\documentclass[reqno,10pt,a4paper]{amsart}
\usepackage{graphicx,color,amssymb,latexsym,amsmath,amsfonts}
\usepackage{mathrsfs}
\numberwithin{equation}{section} \theoremstyle{plain}
\theoremstyle{definition}

\begin{document}

\title[Large deviations for two scaled diffusions]
{Large deviations for two scaled diffusions}

\author{Liptser R.}
\address{Dept. of Electrical Engineering-Systems
Tel Aviv University}
\email{liptser@eng.tau.ac.il}

\date{}

\maketitle
\begin{abstract}
{We formulate large deviations principle (LDP)
for diffusion pair
$(X^\varepsilon,\xi^\varepsilon)=(X_t^\varepsilon,\xi_t^\varepsilon)$,
where first component has a small
diffusion parameter while the second is ergodic Markovian process with
fast time. More exactly, the LDP is established for
$(X^\varepsilon,\nu^\varepsilon)$ with $\nu^\varepsilon(dt,dz)$
being an occupation type measure corresponding to $\xi_t^\varepsilon$. In some
sense we obtain a combination of Freidlin-Wentzell's and Donsker-Varadhan's
results.
Our approach
relies the concept of the exponential tightness and Puhalskii's theorem.}
\end{abstract}

\section{\bf Introduction}

Let $\varepsilon$ be a small positive parameter, $(X^\varepsilon,\xi^
\varepsilon)=(X^\varepsilon_t,\xi_t^\varepsilon)_{t\ge 0}$ be a diffusion
pair defined on some stochastic basis $(\Omega,{\mathcal  F},{\mathsf F}=({\mathcal  F}_t)_
{t\ge 0},{\mathsf P})$ by It\^o's equations w.r.t. independent Wiener
processes $W_t$ and $V_t$:
$$
\begin{aligned}
dX^\varepsilon_t&=A(X_t^\varepsilon,\xi_t^\varepsilon)dt+
\sqrt{\varepsilon}B(X_t^\varepsilon,\xi_t^\varepsilon)dW_t
\\
d\xi_t^\varepsilon&=\frac{1}{\varepsilon}b(\xi_t^\varepsilon)dt+
\frac{1}{\sqrt{\varepsilon}}\sigma(\xi_t^\varepsilon)dV_t
\end{aligned}
\eqno(1.1)
$$
subject to fixed initial point $(x_0,z_0)$.

\medskip
\noindent
Assume $(X^\varepsilon,\xi^\varepsilon)$ is an ergodic process
in the following sense. Let $p(z)$ be the unique invariant density of
$\xi^\varepsilon$,
$$
\nu^{(p)}(dt,dz)=p(z)dtdz,
$$
and $\bar{X}_t$ is a solution of an ordinary
differential equation $\dot{\bar{X}}_t=\bar{A}(\bar{X}_t)$
with $\bar{A}(x)=\int_\mathbb{R}A(x,z)p(z)dz$ subject to the same initial point
$x_0$ . Then for any bounded continuous function $h(t,z)$ and $T>0$
$$
\begin{aligned}
&
{\mathsf P}\text{-}\lim_{\varepsilon\to 0}\int_0^Th(t,\xi_t^\varepsilon)dt=
\int_0^T\int_\mathbb{R}h(t,z)\nu^{(p)}(dt,dz),
\\
&
{\mathsf P}\text{-}\lim_{\varepsilon\to 0}r_T(X^\varepsilon,\bar{X})=0,
\end{aligned}
\eqno(1.2)
$$
where $r_T$ is the uniform metric on $[0,T]$.
The above-mentioned  ergodic property is a motivation to examine
LDP for pair
$(X^\varepsilon,\xi^\varepsilon)$, or more exactly for pair
$(X^\varepsilon,\nu^\varepsilon)$, where $\nu^\varepsilon=
\nu^\varepsilon(dt,dz)$ is an occupation measure on
$\big(\mathbb{R}_+\times \mathbb{R}, {\mathcal  B}(\mathbb{R}_+)\otimes{\mathcal B}
(\mathbb{R})\big)$
(${\mathcal  B}(\mathbb{R}_+)$, ${\mathcal  B}(\mathbb{R})$ are the Borel
$\sigma$-algebras on $\mathbb{R}_+$ and $\mathbb{R}$ respectively) corresponding
to $\xi^\varepsilon$:
$$
\nu^\varepsilon(\Delta\times \Gamma)=\int_0^\infty I(t\in \Delta, \;
\xi_t^\varepsilon\in\Gamma)dt, \quad \Delta\in{\mathcal  B}(\mathbb{R}_+), \;
\Gamma\in{\mathcal  B}(\mathbb{R}).
\eqno (1.3)
$$

\medskip
\noindent
A choice of $\nu^\varepsilon$ as
the occupation measure is natural since the first ergodic property in (1.2)
is nothing but
$$
{\mathsf P}-\lim_{\varepsilon\to 0}\rho_T(\nu^\varepsilon,\nu^{(p)})=0,
$$
where $\rho_T$ is Levy-Prohorov's distance for
restrictions of measures $\nu^\varepsilon$ and $\nu^{(p)}$
on $[0,T]\times \mathbb{R}.$ Also the first It\^o's equation in (1.1)
and the predictable quadratic variation $\big<M^\varepsilon\big>_t$ of
a martingale $M^\varepsilon_t=\int_0^tB(X^\varepsilon_s,\xi^\varepsilon_s)
dW_s$ can be represented in the term of $\nu^\varepsilon$:
$$
\begin{aligned}
&
X_t^\varepsilon=x_0+\int_0^t\int_\mathbb{R}A(X_s^\varepsilon,z)\nu^\varepsilon(ds,dz)+
\sqrt{\varepsilon}M_t^\varepsilon,
\\
&
\big<M^\varepsilon\big>_t
=\int_0^t\int_\mathbb{R}B^2(X_s^\varepsilon,z)\nu^\varepsilon(ds,dz).
\end{aligned}
$$

The random measure $\nu^\varepsilon$
obeys the disintegration
$\nu^\varepsilon(dt,dz)=dtK_{\nu^\varepsilon}(t,dz)$ with the transition kernel
$K_{\nu^\varepsilon}(t,dz)$ being probabilistic Dirac's measure
that is $\nu^\varepsilon$ values
in space $\Bbb{M}=\Bbb{M}_{[0,\infty)}$ of $\sigma$-finite (locally in $t$)
measures $\nu=\nu(dt,dz)$ on $(\mathbb{R}_+\times \mathbb{R}, {\mathcal  B}
(\mathbb{R}_+)\otimes{\mathcal  B}(\mathbb{R}))$
obeying the disintegration $\nu(dt,dz)=K_\nu(t,dz)dt$
with the probabilistic transition kernel $K_\nu(t,dz)$
($\int_\mathbb{R}K_\nu(t,dz)\equiv 1$).
$X^\varepsilon$ values in the space $\Bbb{C}=\Bbb{C}_{[0,\infty)}$
of continuous function.
Define metrics $r$ and $\rho$ in
$\Bbb{C}$ and $\Bbb{M}$ respectively, letting
$$
r(X',X'')=
\sum_{k\ge 1}{ r_k(X',X'')\wedge \frac{1}{2^k}}
\quad\text{and}\quad
\rho(\nu',\nu'')=\sum_{k\ge 1}{
\rho_k(\nu',\nu'')\wedge \frac{1}{2^k}}.
$$
Evidently  ergodic properties
(1.2) are equivalent to
$$
\mathsf{P}-\lim_{\varepsilon\to 0}\big[r(X^\varepsilon,\bar{X})+
\rho(\nu^\varepsilon,\nu^{(p)})\big]=0
$$
and so for examination of the LDP for
$(X^\varepsilon,\nu^\varepsilon)$ we choose the metric space
$(\Bbb{C}\times\Bbb{M},r\times\rho)$.
\medskip
\noindent
Recall the definition of LDP from Varadhan, \cite{1}, adapted to our setting.
The family $(X^\varepsilon,\nu^\varepsilon)$ obeys the LDP
in the metric space $(\Bbb{C}\times\Bbb{M},r\times\rho)$ if

{\bf (0)} there exists a function $L(X,\nu), \;
X\in\Bbb{C},\nu\in\Bbb{M}$,
values in $[0,\infty]$, such that its level sets are compacts in
$(\Bbb{C}\times\Bbb{M},r\times\rho)$;

{\bf (1)}for any open set $G$ from
$(\Bbb{C}\times\Bbb{M},r\times\rho)$
$$
\varliminf_{\varepsilon\to 0}\varepsilon\log {\mathsf P}
\big((X^\varepsilon, \nu^\varepsilon)\in G\big)\ge
-\inf_{(X,\nu)\in G}L(X,\nu);
$$

{\bf (2)} for any closed set $F$ from
$(\Bbb{C}\times\Bbb{M},r\times\rho)$
$$
\varlimsup_{\varepsilon\to 0}\varepsilon\log {\mathsf P}\big((X^\varepsilon,
\nu^\varepsilon)\in F\big)\le -\inf_{(X,\nu)\in F}L(X,\nu).
$$

\noindent
The function $L(X,\nu)$, meeting in {\bf (0)}, {\bf (1)} and {\bf (2)},
is named rate function (action functional in the terminology of Freidlin
and Wentzell, \cite{2}, or good rate function in the terminology of
Stroock, \cite{3}.

\medskip
\noindent
Below we recall well known particular results in LDP's related to pair
$(X^\varepsilon,\xi^\varepsilon)$ and give corresponding forms of rate
functions which will be inherited by a rate function for our setting.
Note at first LDP for family $\mu^\varepsilon(dz)=
\nu^\varepsilon([0,1],dz)$ (on the space of probability measures
supplied by Levy-Prohorov's metric)
 proved by Donsker and Varadhan \cite{4}, \cite
{5}, \cite{6}, \cite{7}
for a wide class of Markov processes $\xi_t^\varepsilon=\xi_{t/\varepsilon}$.
Corresponding rate function obeys an invariant form: for any probabilistic
measure $\mu$ on $\mathbb{R}$
$$
I(\mu)=-\inf \int_\mathbb{R}{{\mathcal  L}u(z)\over u(z)}\mu(dz),
$$
where ${\mathcal  L}$ is backward Kolmogorov's operator, respecting to
$\xi$, and where `$inf$' is taken over all functions $u(z)$ from the domain of
definition for the operator $\mathcal {L}$. For the
diffusion case, G\"artner's
type of $I(\mu)$ is well known (\cite{8}):
$$
I(\mu)=
\begin{cases}

{1\over 8}\int_\mathbb{R}\sigma^2(z)\Big[{m'(z)\over m(z)}-{p'(z)\over p(z)}\Big]^2m(z)dz,
& d\mu=m(z)dz, \ dm(z)=m'(z)dz,\\
\infty, & \text{otherwise}.
\end{cases}
\eqno (1.4)
$$
Freidlin-Wentzell's result, \cite{2}, is devoted
to LDP for diffusion $X^\varepsilon$ with drift $A(x)$ and
diffusion $B^2(x)$ (independent of $z$) in the space of continuous functions
on every finite time interval, supplied by the uniform metric.
A rate function, say, for $[0,T]$ time interval is given by
$$
S(X)=
\begin{cases}
{1\over 2}\int_0^T{[\dot{X}_t-A(X_t)]^2\over B^2(X_t)}dt, &
dX_t=\dot{X}_tdt, \ X_0=x_0\\
\infty, & \text{otherwise}.
\end{cases}
\eqno (1.5)
$$
Other type of LDP for a degenerate diffusion $X^\varepsilon$ defined
by the first equation in (1.1) with $B(x,z)\equiv 0$ and $\xi_t^\varepsilon=
\xi_{t/\varepsilon}$, where $\xi_t$ is Markov process values in a finite
state space, also is well known from Freidlin \cite{9}. In this case
rate function has
a form similar to (1.5) ($H(y,x)$ is some non negative function):
$$
S(X)=
\begin{cases}
\int_0^TH(\dot{X}_t,X_t)dt, &
dX_t=\dot{X}_tdt, \ X_0=x_0\\
\infty, & \text{otherwise}.
\end{cases}
\eqno (1.6)
$$

\medskip
\noindent
All above-mentioned LDP's  are inspired the examination of the LDP
for $(X^\varepsilon,\nu^\varepsilon)$.
In some sense, the LDP for $(X^\varepsilon,\nu^\varepsilon)$ is a
combination of Donsker-Varadhan's and Freidlin-Wenzell's results. Namely
LDP for $\nu^\varepsilon$ is a generalization one for
$\mu^\varepsilon$ while LDP for $X^\varepsilon$ is implied by
LDP for $\nu^\varepsilon$ and for a diffusion martingale scaled by
$\sqrt{\varepsilon}$. Hence, a rate function for
$(X^\varepsilon,\nu^\varepsilon)$, is defined as a sum: $L(X,\nu)=
L_1(X,\nu)+L_2(\nu),$ where $L_1(X,\nu)$ and $L_2(\nu)$ respect
to $X^\varepsilon$ and $\nu^\varepsilon$
and what is more  $L_1(X,\nu)$ has the same form as
$S(x)$ in (1.5) with $A(X_t)$ and $B^2(X_t)$ replaced on
$A_\nu(t,X_t)=\int_\mathbb{R}A(X_t,z)K_\nu(t,dz)$ and
$B^2_\nu(t,X_t)=\int_\mathbb{R}B^2(X_t,z)K_\nu(t,dz)$, where $K_\nu(t,dz)$
is the transition kernel of measure $\nu.$

\medskip
\noindent
Note that $\xi^\varepsilon_t\in R$ and so
the LDP for its occupation measure responds to a non compact diffusion case.
Also note that diffusion parameter $B^2(x,z)$ is not assumed to be non
singular and
consequently $B^2(x,z)\equiv 0$ is admissible. The last allows to derive
LDP for a singular diffusion parameter case
from LDP for $\nu^\varepsilon$ using the contraction principle
of Varadhan \cite{1} (continuous mapping method of Freidlin \cite{10}).
This result extends above-mentioned \cite{9} for
non compact case.

\medskip
\noindent
In contrast with Freidlin and Wentzell \cite{2}, Donsker and Varadhan
\cite{4} -\cite{7}, G\"artner \cite{8}, and Veretennikov \cite{11} - \cite{12},
and many
others (see e.g. Acosta \cite{13}, Dupuis and Elis \cite{14}) our method
of proof is based on Puhalskii's theorem \cite{15} - \cite{16} and relies concepts
of exponential tightness and LD relative compactness.

\medskip
\noindent
The paper is organized as follows. In Section 2, we formulate
the general assumptions and the main result. Section 3 contains the
method of proving LDP which also has been used in \cite{17}.
In Section 4, we check the exponential tightness while in Section 5 and 6
the upper and lower bounds in local LDP are verified. The main
results are proved in Section 7.
All technical results are gathered in Appendix.

\section{\bf Assumptions. Main result}
\mbox{}

\noindent
{\bf 1.} We fix the following conditions which are assumed to be fulfilled
hereafter.

{\bf (A.1)} $A(x,z)$ and $B(x,z)$ are continuous in $(x,z)$,
Lipschitz continuous in $x$ uniformly in $z$, and $\sup_z\big(|A(0,z)|+
|B(0,z)|\big)<\infty$;

{\bf (A.2)} $\sigma^2(z)$ is bounded and and uniformly positive
function; it is continuously
differentiable, having bounded and Lipschitz continuous derivative;

{\bf (A.3)} $b(z)$ is Lipschitz continuous, satisfying
$$
\lim_{|z|\to\infty}b(z)\text{sign} \ z=-\infty.
$$

\medskip
\noindent
It would be noted that (A.2) and (A.3) imply, so called, assumption
(H$^*$) from \cite{6}.

\medskip
\noindent
{\bf 2.} It is well known (see \cite{18}) that under (A.2) and (A.3)
$\xi^\varepsilon$ is ergodic process obeying the unique invariant density
$$
p(z)=\text{const.}{\exp\Big(2\int_0^z\frac{b(y)}{\sigma^2(y)}dy\Big)
\over\sigma^2(z)}.
\eqno (2.1)
$$
For any $\nu$ from $\Bbb{M}$ with the transition kernel $K_\nu(t,dz)$, define
$K_\nu(t,dz)$-averaged drift $A_\nu(t,x)$ $=\int_\mathbb{R}A(x,z)K_\nu(t,dz)$ and
diffusion parameter
$$
B^2_\nu(t,x)=\int_\mathbb{R}B^2(x,z)K_\nu(t,dz).
$$
If $\nu$ is absolutely continuous w.r.t.
$\Lambda(dt,dz)=dtdz$,  put
$$
n(t,z)={d\nu\over d\Lambda}(t,z).
\eqno (2.2)
$$
If the density $n(t,z)$ is absolutely continuous w.r.t. $dz$:
$d_zn(t,z)=n'_z(t,z)dz$, a function $n'_z(t,z)$ is chosen to be measurable
in $t,z.$

\medskip
\noindent
Throughout the paper, we use conventions $0/0=0$ and $\min(\inf)(\emptyset)=
\infty$.

\medskip
\noindent
For every $\nu\in \Bbb{M}$ and $X\in\Bbb{C}$
define two quantities (comp. (1.4) and (1.5)):
$$
\begin{aligned}
F(\nu)&=
\begin{cases}
\int_0^\infty\int_\mathbb{R}\sigma^2(z)\Big[{n'_z(t,z)\over n(t,z)}-
{p'(z)\over p(z)}\Big]^2n(t,z)dzdt,
& d\nu=nd\Lambda, \ d_zn=n'_zdz
\\
\infty, & \text{otherwise};
\end{cases}
\\
S(X,\nu)&=
\begin{cases}
\int_0^\infty{[\dot{X}_t-A_\nu(t,X_t)]^2\over B^2_\nu(t,X_t)}dt, &
dX=\dot{X}dt, \ \ \ X_0=x_0
\\
\infty, & \text{otherwise}.
\end{cases}
\end{aligned}
\eqno(2.3)
$$

\medskip
\noindent
{\bf 3.} Now we are in the position to formulate the main result.

{\bf Theorem 2.1.} {\it Under {\rm (A.1)}, {\rm (A.2)}, and {\rm (A.3)} the family
$(X^\varepsilon,\nu^\varepsilon)$ obeys the LDP in
$(\Bbb{C}\times\Bbb{M},r\times\rho)$ with rate function
$$
L(X,\nu)={1\over 2}S(X,\nu)+{1\over 8}F(\nu).
$$
}

\medskip
\noindent
{\bf 4.} LDP's for families $(X^\varepsilon)$ and $(\xi^\varepsilon)$
run out from Theorem 2.1.

{\bf Corollary 2.1.} {\it $(\nu^\varepsilon)$ obeys the LDP in
$(\Bbb{M},\rho)$ with rate function ${1\over 8}F(\nu).$}

\medskip
{\bf Corollary 2.2.} {\rm (comp. \cite{9})}{\it
$(X^\varepsilon)$ obeys the LDP in
$(\Bbb{C},r)$ with rate function $S(X)=\inf_{\nu\in\Bbb{M}}L(X,\nu).$
In particular, if $B(x,z)\equiv 0$, it is sufficiently to take
`$inf$' over all $\nu$ from $\Bbb{M}$ with the transition kernel
$K_\nu(t,dz)\equiv \mu(dz)$
with $d\mu=m(z)dz$ such that
the density $m(z)={d\mu\over dz}(z)$ is absolutely continuous
w.r.t. $dz$ ($m'(z)={dm(z)\over dz}$). In this case, rate function
$$
S(X)=
\begin{cases}
{1\over 8}\int_0^\infty H(\dot{X}_t,X_t)dt, & dX=\dot{X}dt, \ \ \ X_0=x_0
\\
\infty, & \text{otherwise},
\end{cases}
\eqno (2.4)
$$
where
$$
H(y,x)=\inf\int_\mathbb{R}\sigma^2(z)\Big[{m'(z)\over m(z)}
-{p'(z)\over p(z)}\Big]^2m(z)dz,
\eqno (2.5)
$$
and where `$inf$' is taken over all above-mentioned measures
$\mu$ such that
$$
y=\int_\mathbb{R}A(x,z)m(z)dz.
$$}

\medskip
\noindent
As an example, also the LDP for the family of the Donsker and Varadhan
occupation measures $\mu^\varepsilon(dz)=\nu^\varepsilon([0,1]\times dz)$,
corresponding to diffusion case,
can be derived from Theorem 2.1. In fact,
due to the contraction principle, $(\mu^\varepsilon)$ obeys the LDP with
G\"artner's type rate function (see (1.4)) $I(\mu)=\inf{1\over 8}F(\nu),$
where `$inf$' is taken over all $\nu\in\Bbb{M}$ such that
$$
\nu(dt,dz)=I(1\ge t)dt\mu(dz)+I(1<t)\nu^{(p)}(dt,dz).
$$

\section{\bf Preliminaries}

For proving LDP for the family $(X^\varepsilon,\nu^\varepsilon)$
in the metric space $(\Bbb{C}\times\Bbb{M},r\times\rho)$
we apply Dawson-G\"artner's type theorem (see e.g. \cite{19}. Following it
the LDP in $(\Bbb{C}\times\Bbb{M},r\times\rho)$
is implied by LDP's in the metric spaces
$(\Bbb{C}_{[0,n]}\times\Bbb{M}_{[0,n]},r_n\times\rho_n), \; n\ge 1$, where
$\Bbb{C}_{[0,n]}$ is the space of continuous functions on the time interval
$[0,n]$, $\Bbb{M}_{[0,n]}$ is the space of finite measures on $[0,n]\times \mathbb{R}$,
having probabilistic transition kernel w.r.t. $dt$, $r_n$ is the uniform
metric, and $\rho_n$ is Levy-Prohorov's metric. The definition
of the LDP in $(\Bbb{C}_{[0,n]}\times\Bbb{M}_{[0,n]},r_n\times\rho_n)$
is given in terms of ({\bf 0.}), ({\bf 1.}), and ({\bf 2.})
with obvious modifications. Moreover, if $L_n(X,\nu), n\ge 1$ are
rate functions, corresponding to LDP's in
$(\Bbb{C}_{[0,n]}\times\Bbb{M}_{[0,n]},r_n\times\rho_n), n\ge 1$, then rate
function in $(\Bbb{C}\times\Bbb{M},r\times\rho)$ is defined as
$$
L(X,\nu)=\sup_nL_n(X,\nu).
\eqno (3.1)
$$
Hence only the LDP in
$(\Bbb{C}_{[0,T]}\times\Bbb{M}_{[0,T]},r_T\times\rho_T)$ has to be checked
for any $T>0$. Our approach in proving the LDP in
$(\Bbb{C}_{[0,T]}\times\Bbb{M}_{[0,T]},
r_T\times\rho_T), T>0$ relies the concept of the exponential tightness
and notions of LD relative compactness and local LDP. Below we give necessary
definitions.

\medskip
\noindent
{\bf Definition 1.} The family
$(X^\varepsilon,\nu^\varepsilon)$ is said to be exponentially tight
in the metric space
$(\Bbb{C}_{[0,T]}\times\Bbb{M}_{[0,T]},r_T\times\rho_T)$, if there
exists an increasing sequence of compacts $(K_j)_{j\ge 1}$ such that
$$
\lim_j\bar{\lim_{\varepsilon\to 0}}\varepsilon\log
{\mathsf P}\big((X^\varepsilon,\nu^\varepsilon)\in
\{\Bbb{C}_{[0,T]}\times\Bbb{M}_{[0,T]}\}\setminus K_j\big)=-\infty.
\eqno (3.2)
$$
\newline (Deuschel and Stroock \cite{20}, Lynch and Sethuraman \cite{21}.)

\medskip
\noindent
{\bf Definition 2.} The family
$(X^\varepsilon,\nu^\varepsilon)$ is said to be
LD relatively compact in
$(\Bbb{C}_{[0,T]}\times\Bbb{M}_{[0,T]},r_T\times\rho_T)$, if any
decreasing to zero sequence $(\varepsilon_k)$ contains further subsequence
$(\bar{\varepsilon}_k)$ ($(\bar{\varepsilon}_k)\subseteq
(\varepsilon_k)$) such that the family
$(X^{\bar{\varepsilon}_k},\nu^{\bar{\varepsilon}_k})$
obeys the LDP in
$(\Bbb{C}_{[0,T]}\times\Bbb{M}_{[0,T]},r_T\times\rho_T)$
(with rate function $\bar{L}_T(X,\nu)$).
\newline (Puhalskii \cite{15} - \cite{16}.)

\medskip
\noindent
{\bf Definition 3.}
The family $(X^\varepsilon,\nu^\varepsilon)$ is said to be obey
the local LDP in
$(\Bbb{C}_{[0,T]}\times\Bbb{M}_{[0,T]},r_T\times\rho_T)$
with local
rate function $\widehat{L}_T(X,\nu)$, if for any $(X,\nu)$ from
$\Bbb{C}_{[0,T]}\times\Bbb{M}_{[0,T]}$
$$
\begin{aligned}
&\quad\varlimsup_{\delta\to 0}
\varlimsup_{\varepsilon\to 0}\varepsilon\log
{\mathsf P}\big((r_T(X^\varepsilon,X)+\rho_T(\nu^\varepsilon,\nu)\le
\delta\big)
\\
&
=\varliminf_{\delta\to 0}\varliminf_{\varepsilon\to 0}\varepsilon\log
{\mathsf P}\big((r_T(X^\varepsilon,X)+\rho_T(\nu^\varepsilon,\nu)\le
\delta\big)
\\
&
=-\widehat{L}_T(X,\nu).
\end{aligned}
\eqno(3.3)
$$
\newline (Freidlin and Wentzell \cite{2}.)

\medskip
\noindent
The connecting component of these notions used in the proof  of the
next result is Puhalskii's theorem \cite{15}, \cite{16}. Below we formulate only
the first part of it.

{\bf Theorem P.} {\it If $(X^\varepsilon,\nu^\varepsilon)$
is exponentially tight family in
$(\Bbb{C}_{[0,T]}\times\Bbb{M}_{[0,T]},r_T\times\rho_T)$,
then it is LD relatively compact.}

\medskip
\noindent
The following result is a reformulation of Theorem 1.3 from \cite{17}.

\medskip
{\bf Proposition 3.1.} {\it The exponential tightness and the local LDP
for the family
$(X^\varepsilon,\nu^\varepsilon)$
in $(\Bbb{C}_{[0,T]}\times\Bbb{M}_{[0,T]},r_T\times\rho_T)$
imply the LDP in
$(\Bbb{C}_{[0,T]}\times\Bbb{M}_{[0,T]},r_T\times\rho_T)$
for this family with (good) rate function
$L_T(X,\nu)\equiv \widehat{L}_T(X,\nu)$, where $\widehat{L}_T(X,\nu)$
is the local rate function.}

\section{\bf Exponential tightness in $\boldsymbol{\Bbb{C}_{[0,T]}\times
\Bbb{M}_{[0,T]}}$}
\mbox{}

{\bf Theorem 4.1.} {\it Under assumptions (A.1), (A.2), and (A.3)
the family $(X^\varepsilon,\nu^\varepsilon)$ is
exponentially tight in $\Bbb{C}_{[0,T]}\times\Bbb{M}_{[0,T]}$.}

\medskip
{\it Proof.} Following Definition 1, (3.2) has to be checked. It is clear
it takes place if
$$
\begin{aligned}
&
\lim_j\varlimsup_{\varepsilon\to 0}\varepsilon\log {\mathsf P}\big(
X^\varepsilon\in \Bbb{C}_{[0,T]}\setminus K_j'\big)=-\infty
\\
&
\lim_j\varlimsup_{\varepsilon\to 0}\varepsilon\log {\mathsf P}\big(
\nu^\varepsilon\in \Bbb{M}_{[0,T]}\setminus K_j''\big)
=-\infty,
\end{aligned}
\eqno(4.1)
$$
where $K_j'$ and $K_j''$ are appropriate increasing sequences of compacts
from $\Bbb{C}_{[0,T]}$ and $\Bbb{M}_{[0,T]}$ respectively.
It is naturally to use as compacts $K_j'$
increasing sets of uniformly bounded and equicontinuous functions
from $\Bbb{C}_{[0,T]}$ parametrized by $j$. Since the process
$(X^\varepsilon_t, \xi^\varepsilon_t)_{t\ge 0}$
is defined on a stochastic basis with the filtration ${\mathsf F}$
one can use Aldous-Puhalskii's
type sufficient conditions
(see \cite{15}, and also Theorem 3.1 in
\cite{17}) for $C$-exponential tightness:
$$
\begin{aligned}
&
\lim_j\bar{\lim_{\varepsilon\to 0}}\varepsilon\log {\mathsf P}\big(
\sup_{t\le T}|X^\varepsilon_t|>j\big)=-\infty
\\
&
\lim_{\delta\to 0}\bar{\lim_{\varepsilon\to 0}}\varepsilon\log
\sup_{\tau\le T-\delta}{\mathsf P}\big(\sup_{t\le \delta}|X^\varepsilon_{\tau+t}-
X_\tau^\varepsilon|>\eta\big)=-\infty, \; \forall \eta>0,
\end{aligned}
\eqno(4.2)
$$
where $\tau$ is a stopping time w.r.t. the filtration ${\mathsf F}$. Following
Theorem 3.1 in \cite{17}, (4.2) implies the validity of the first part
in (4.1) with above-mentioned compacts $K'_j$ of uniformly bounded and
equicontinuous functions. Now, choose relevant compacts $K''_j, j\ge 1$:
$$
K_j''=\bigcap_{m\ge j}\big\{\nu\in\Bbb{M}_{[0,T]}:\int_0^T\int_{|z|>m}
\nu(dt,dz)\le g(m)\big\},
\eqno (4.3)
$$
where $g(y), y>0$ is positive continuous decreasing function with
$\lim_{y\to\infty}g(y)=0$.
In fact, if $\nu_k\in K_j'',k\ge 1$ then we have for any $m\ge j$,
$\sup_k\int_0^T\int_{|z|>m}\nu_k(dt,dz)\le g(m)$ that is the set
$K_j''$ is tight and by Prohorov's theorem (see \cite{22}) is relatively
compact. On the other hand, since the set $\{z:|z|>m\}$ is open a limit of any
converging sequence from $K_j''$ also belongs to $K_j''$ that is
$K_j''$ is compact
in $(\Bbb{M}_{[0,T]},\rho_T)$. Evidently $K_j''\subseteq K_{j+1}''$.
Below we choose a special function $g(y)$, suited to assumption (A.3),
to satisfy the second part in (4.1).
\medskip
\noindent
We check the
validity of (4.1) in the next two lemmas.

\smallskip
{\bf Lemma 4.1.} {\it Under (A.1) the first relation in (4.1) holds.}

\smallskip
{\bf Lemma 4.2.} {\it Under (A.2) and (A.3) the second relation in (4.1)
holds.}

\medskip
{\it Proof of Lemma 4.1.} Put
$$
Z^*_t=\sup_{t'\le t}|Z_{t'}|.
$$

\medskip
\noindent
By virtue of (A.1) we have $|A(x,z)|\le \ell(1+|x|)$. Therefore,
with $t\le T$, we derive from (1.1)
$$
X_t^{\varepsilon *}\le |x_0|+\ell\int_0^t(1+X_s^{\varepsilon *})ds+
\sqrt{\varepsilon}M_T^{\varepsilon *},
\eqno (4.4)
$$
where $M_t^\varepsilon=\int_0^tB(X_s^\varepsilon,
\xi_s^\varepsilon)dW_s.$ Due to Bellman-Gronwall's inequality,
(4.4) implies $X_T^{\varepsilon *}\le\text{const.}(1+\sqrt{\varepsilon}
M_T^{\varepsilon *})$
with const., depending only on $|x_0|, \ell,$ and $T$. Therefore, the first
part of (4.2) holds if
$$
\lim_j\bar{\lim_{\varepsilon\to 0}}\varepsilon\log {\mathsf P}\big(
M_T^{\varepsilon *}>j\big)=-\infty.
\eqno (4.5)
$$
On the other hand, by Chebyshev's inequality ${\mathsf P}\big(
M_T^{\varepsilon *}>j\big)\le j^{-1/\varepsilon}{\mathsf E}\big(
M_T^{\varepsilon *}\big)^{1/\varepsilon}$ and so,
$\varepsilon\log {\mathsf P}\big(M_T^{\varepsilon *}>j\big)\le
-\log j+
\varepsilon\log {\mathsf E}\big(M_T^{\varepsilon *}\big)^{1/\varepsilon}.$
Thereby (4.5) holds if
$$
\bar{\lim_{\varepsilon\to 0}}\varepsilon \log{\mathsf E}
\big(M_T^{\varepsilon *}\big)^{1/\varepsilon}< \infty.
\eqno (4.6)
$$

\medskip
\noindent
Below we check the validity (4.6). Assuming $1/\varepsilon >2$ and
applying It\^o's formula to
$|M_t^\varepsilon|^{1/\varepsilon}$, we get
$$
\begin{aligned}
|M_t^\varepsilon|^{1/\varepsilon}
&
={1\over\sqrt{\varepsilon}}\int_0^t|M_s^\varepsilon|^
{1/\varepsilon-1} \ (\text{sign} \ M_s^\varepsilon) B(X_s^\varepsilon,
\xi_s^\varepsilon)dW_s
\\
&
+{1-\varepsilon\over 2\varepsilon}
\int_0^t|M_s^\varepsilon|^{1/\varepsilon-2}B^2(X_s^\varepsilon,
\xi_s^\varepsilon)ds
\end{aligned}
$$
that is $|M_t^\varepsilon|^{1/\varepsilon}$ is a submartingale obeying
a decomposition: $|M_t^\varepsilon|^{1/\varepsilon}
=N_t^\varepsilon + U_t^\varepsilon$ with local martingale
$N_t^\varepsilon$ and predictable increasing process
$$
U_t^\varepsilon={1-\varepsilon\over 2\varepsilon}
\int_0^t|M_s^\varepsilon|^{1/\varepsilon-2}B^2(X_s^\varepsilon,
\xi_s^\varepsilon)ds.
\eqno (4.7)
$$
Then, due to a modification of Doob's inequality (see \cite{23},
Theorem 1.9.2)
$$
{\mathsf E}\big(M_t^{\varepsilon *}\big)^{1/\varepsilon}\le
\Big({1\over 1-\varepsilon}\Big)^{1/\varepsilon}
{\mathsf E}U^\varepsilon_t.
\eqno (4.8)
$$
Now evaluate from above $|M_s^\varepsilon|^{1/\varepsilon-2}B^2
(X_s^\varepsilon,\xi_s^\varepsilon)$.
By virtue of (A.1)
$|B(x,z)|\le \ell(1+|x|)$. Thereby, due to
above-mentioned upper bound
$X_T^{\varepsilon *}\le \text{const.}(1+M_T^{\varepsilon *})$
which remains true with replacing $T$ on $s$ for any $s<T$,
we arrive at
$$
\begin{aligned}
|M_s^\varepsilon|^{1/\varepsilon-2}B^2(X_s^\varepsilon,
\xi_s^\varepsilon)
&
\le \text{const.}(1+|M_s^\varepsilon|^{1/\varepsilon-2}+
|M_s^\varepsilon|^{1/\varepsilon})
\\
&
\le \text{const.}(1+(M_s^{\varepsilon *})^{1/\varepsilon}).
\end{aligned}
$$
Substituting  the last upper bound in (4.7) and using (4.8) we obtain
($t\le T$)
\newline
${\mathsf E}(M_t^{\varepsilon *})^{1/\varepsilon}\le
{\text{const.}\over \varepsilon}
\int_0^t\big[1+{\mathsf E}(M_t^{\varepsilon *})^{1/\varepsilon}\big]ds$.
Hence, by Bellman-Gronwall's inequality, an upper bound
${\mathsf E}(M_T^{\varepsilon *})^{1/\varepsilon}\le{\text{const.}T
\over \varepsilon}\exp\{{\text{const.}T\over \varepsilon}\}$ holds and
implies (4.6). Consequently the first part in (4.2) is valid.
To check the second part in (4.2), first use obvious estimates:
\begin{multline*}
{\mathsf P}\big(\sup_{t\le \delta}|X^\varepsilon_{\tau+t}-
X_\tau^\varepsilon|>\eta\big)
\le {\mathsf P}\big(\sup_{t\le \delta}|X^\varepsilon_{\tau+t}-
X_\tau^\varepsilon|>\eta, X^{\varepsilon *}_T\le j\big)+
{\mathsf P}\big(X^{\varepsilon *}_T>j\big)
\\
\le 2\max\Big[
{\mathsf P}\big(\sup_{t\le \delta}|X^\varepsilon_{\tau+t}-
X_\tau^\varepsilon|>\eta, X^{\varepsilon *}_T\le j\big),
{\mathsf P}\big(X^{\varepsilon *}_T>j\big)\Big].
\end{multline*}
Thence, due to proved above the first part of (4.2), the validity of the
second part follows if
$$
\lim_{\delta\to 0}\bar{\lim_{\varepsilon\to 0}}\varepsilon\log
\sup_{\tau\le T-\delta}{\mathsf P}\big(\sup_{t\le \delta}|X^\varepsilon_{\tau+t}-
X_\tau^\varepsilon|>\eta, X^{\varepsilon *}_T\le j\big)=-\infty, \ j\ge 1,
\eta >0.
\eqno (4.9)
$$
The simplest way for verifying of (4.9) consists in checking the validity
of both
$$
\begin{aligned}
&
\lim_{\delta\to 0}\bar{\lim_{\varepsilon\to 0}}\varepsilon\log
\sup_{\tau\le T-\delta}{\mathsf P}\big(\sup_{t\le \delta}
\Big|\int_\tau^{\tau+t}A(X_s^\varepsilon,\xi_s^\varepsilon)ds\Big|>\eta,
X^{\varepsilon *}_T\le j\big)=-\infty
\\
&
\lim_{\delta\to 0}\bar{\lim_{\varepsilon\to 0}}\varepsilon\log
\sup_{\tau\le T-\delta}{\mathsf P}\big(\sup_{t\le \delta}
\Big|\sqrt{\varepsilon}\int_\tau^{\tau+t}
B(X_s^\varepsilon,\xi_s^\varepsilon)dW_s\Big|>\eta,
X^{\varepsilon *}_T\le j\big)=-\infty.
\end{aligned}
\eqno(4.10)
$$
Obviously, the first part in (4.10)  holds. To verify the second, note that
the process  $Y_t^\varepsilon=\sqrt{\varepsilon}\int_\tau^{\tau+t}
B(X_s^\varepsilon,\xi_s^\varepsilon)dW_s$ is continuous martingale
w.r.t. the new filtration ${\mathsf F}^\tau=(\mathcal {F}_{\tau+t})_{t\ge 0}$
(see Ch.4, \S 7 in \cite{23}). It has the predictable quadratic variation
$\big<Y^\varepsilon\big>_t=\varepsilon\int_\tau^{\tau+t}
B^2(X_s^\varepsilon,\xi_s^\varepsilon)ds.$
Also define a positive continuous local martingale (w.r.t. the same
filtration ${\mathsf F}^\tau$)
$$
Z_t^\varepsilon=\exp\big(\lambda Y_t^\varepsilon-{1\over 2}\lambda^2
\big<Y^\varepsilon\big>_t\big), \ \lambda\in \mathbb{R}
\eqno (4.11)
$$
which is simultaneously a supermartingale (see \cite{23}, Problem 1.4.4)
and so for any Markov time $\sigma$ (w.r.t. ${\mathsf F}^\tau$)
$\mathsf{E}Z_\sigma^\varepsilon\le 1$. Take
$\sigma=\inf\{t\le\delta:|Y_t^\varepsilon|\ge \eta\}$. Evidently the
second part of (4.10) holds if
$$
\lim_{\delta\to 0}\bar{\lim_{\varepsilon\to 0}}\varepsilon\log
\sup_{\tau\le T-\delta}{\mathsf P}\big(Y_\sigma^\varepsilon\ge \eta
 \ (\text{or} \ \le -\eta),
\sigma\le \delta,
X^{\varepsilon *}_T\le j\big)=-\infty,
\eqno (4.12)
$$
By virtue of an obvious inequality $\mathsf{E}Z_\sigma^\varepsilon
I\big(Y_\sigma^\varepsilon\ge \eta,X^{\varepsilon *}_T\le j\big)\le 1$
we find that
$$
\varepsilon\log{\mathsf P}\big(Y_\sigma^\varepsilon\ge \eta,
\sigma\le \delta,X^{\varepsilon *}_T\le j\big)\le -\sup_{\lambda>0}
\big[\lambda\eta-\text{const.}{\lambda^2\over 2}\delta\varepsilon\big]
\eqno (4.13)
$$
and since $\sup_{\lambda>0}[\lambda\eta-\text{const.}{\lambda^2\over 2}
\delta\varepsilon]={\eta^2\over 2\text{const.}\delta\varepsilon}$ (4.12)
with `$\ge \eta$' is implied by (4.13). The validity (4.12) with
`$\le -\eta$' is proved in the same way.

\medskip
\noindent
{\it Proof of Lemma 4.2.} It is clear that
$\{\nu^\varepsilon\in \Bbb{M}_{[0,T]}\setminus
K_j''\}=\{\ell (j,\nu^\varepsilon)<\infty\}$, where
$$
\ell(j,\nu)=\min\{m\ge j:\int_0^T\int_{|z|>m}\nu(dt,dz)>g(m)\}.
\eqno (4.14)
$$
Therefore, the second part of (4.1) is equivalent to
$$
\lim_j\bar{\lim_{\varepsilon\to 0}}\varepsilon\log {\mathsf P}
(\ell (j,\nu^\varepsilon)<\infty)=-\infty.
\eqno (4.15)
$$
To verify (4.15), choose a special function $g(y)$ satisfying
above-mentioned properties. To this end introduce non linear operator
$$
\mathcal {D}=b(z){\partial\over \partial z}+{\sigma^2(z)\over 2}\big[
{\partial^2\over \partial z^2}+\big({\partial\over \partial z}\big)^2\big]
\eqno (4.16)
$$
and choose a non negative twice continuously differentiable function $u(z)$
such that
$$
\begin{aligned}
&
-\sup_{v\in \mathbb{R}}{\mathcal  D}u(v)=-d>-\infty,
\\
&
\lim_{j\to\infty}\inf_{|z|>j}\big[-{\mathcal  D}u(z)+\sup_{v\in \mathbb{R}}{\mathcal  D}u(v)\big]
=\infty.
\end{aligned}
\eqno(4.17)
$$
Under assumptions (A.2) and (A.3) one can take any of function $u(z)$
with properties:
$u(0)=0$, $u'(z)=\text{sign} \ z, \; |z|>1$, and $0\le u''(z)\le 1$.
With chosen $u(z)$ put
$$
g(y)=\inf_{|z|>y}\big[-{\mathcal  D}u(z)+\sup_v{\mathcal  D}u(v)\big]^{-1/2}.
\eqno (4.18)
$$
Introduce
a positive continuous local martingale
(the martingale property is checked by It\^o's formula)
$$
Z_t^\varepsilon=\exp\big(u(\xi_t^\varepsilon)-u(\xi_0)-
\frac{1}{\varepsilon}\int_0^t{\mathcal  D}u(\xi_s^
\varepsilon)ds\big).
\eqno (4.19)
$$
It is simultaneously a supermartingale (see Problem 1.4.4. in \cite{23})
and so $\mathsf{E}Z_T^\varepsilon\le 1$. The last implies
$$
{\mathsf E}I(\ell(j,\nu^\varepsilon)<\infty)Z_T^\varepsilon\le 1.
\eqno (4.20)
$$
Inequality (4.20) can be sharpen by changing of $Z_T^\varepsilon$
on its lower bound on the set
\newline
$\{\ell(j,\nu^\varepsilon)<\infty\}$
which can be chosen non random.
Taking into account that $\int_0^T{\mathcal  D}u(\xi_s^\varepsilon)ds=
\int_0^T\int_\mathbb{R}{\mathcal  D}u(z)\nu^\varepsilon(ds,dz)$ and
$\ell(j,\nu^\varepsilon)\ge j$, write
$$
\begin{aligned}
\log Z_T^\varepsilon
&
\ge -u(\xi_0)-{dT\over\varepsilon}+
{1\over\varepsilon}\int_0^T\int_{|z|>\ell(j,\nu^\varepsilon)}
\big[-{\mathcal  D}u(z)+d\big]\nu^\varepsilon(ds,dz)
\\
&
\ge -u(\xi_0)-{dT\over\varepsilon}+
{1\over\varepsilon}\inf_{|z|>\ell(j,\nu^\varepsilon)}
\big[-{\mathcal  D}u(z)+d\big]
\int_0^T\int_{|z|>\ell(j,\nu^\varepsilon)}\nu^\varepsilon(ds,dz)
&
\cr
&
\ge -u(\xi_0)-{dT\over\varepsilon}+
{1\over\varepsilon}\inf_{|z|>\ell(j,\nu^\varepsilon)}
\big[-{\mathcal  D}u(z)+d\big]^{1/2}
&
\cr
&
\ge -u(\xi_0)-{dT\over\varepsilon}+
{1\over\varepsilon}\inf_{|z|>j}
\big[-{\mathcal  D}u(z)+d\big]^{1/2} \ (=\log Z_*).
\end{aligned}
$$
Thereby, from (4.20), with $Z_T^\varepsilon$ replaced on $Z_*$, we derive
$$
\varepsilon\log{\mathsf P}\big(\ell(j,\nu^\varepsilon)<\infty)\le
\varepsilon u(\xi_0)+dT-\inf_{|z|>j}\big[-{\mathcal  D}u(z)+d\big]^{1/2},
$$
i.e. (4.15) is implied by (4.17).

\section{\bf Upper bound for local LDP in $\boldsymbol{\Bbb{C}_{[0,T]}\times
\Bbb{M}_{[0,T]}}$}

In this Section, we consider family $(X^\varepsilon,
\nu^\varepsilon)$ from $\Bbb{C}_{[0,T]}\times\Bbb{M}_{[0,T]}$.
Parallel to $F(\nu)$ and
$S(X,\nu)$, given in (2.3), let us define $F_T(\nu)$ and $S_T(X,\nu)$
by changing integrals `$\int_0^\infty$' in (2.3) on
`$\int_0^T$'. Put
$$
L_T(X,\nu)={1\over 2}S_T(X,\nu)+{1\over 8}F_T(\nu).
\eqno (5.1)
$$

\medskip
{\bf Theorem 5.1.} {\it Under {\rm (A.1), (A.2)}, and {\rm (A.3)} for every $(X,\nu)$
from $\Bbb{C}_{[0,T]}\times\Bbb{M}_{[0,T]}$
$$
\bar{\lim_{\delta\to 0}} \ \bar{\lim_{\varepsilon\to 0}}
\varepsilon\log {\mathsf P}\big(r_T(X^\varepsilon,X)+\rho_T(\nu^\varepsilon,
\nu)\le \delta\big)\le -L_T(X,\nu).
$$}

\medskip
Proof of this theorem is based on

\medskip
{\bf Lemma 5.1.} {\it Assume (A.1), (A.2), and (A.3). Then for every
piece wise constant function $\lambda(t)=\sum_i\lambda(t_i)I(t_i\le t
<t_{i+1})$ (with not overlapping intervals $[t_i,t_{i+1})$),
and for every compactly supported in $z$ and continuously
differentiable
(once in $t$ and twice in $z$) function $u(t,z)$, and $X\in \Bbb{C}_{[0,T]},
\ \nu\in \Bbb{M}_{[0,T]}$
$$
\begin{aligned}
\bar{\lim_{\delta\to 0}}
&
 \ \bar{\lim_{\varepsilon\to 0}}
\varepsilon\log{\mathsf P}\big(r_T(X^\varepsilon,X)+\rho_T
(\nu^\varepsilon,\nu)\le\delta\big)
\\
&
\le -\Big\{\sum_i\lambda(t_i)[X_{T\wedge t_{i+1}}-X_{T\wedge t_i}]-
\int_0^T\int_\mathbb{R}\lambda(t)A(X_t,z)\nu(dt,dz)
&
\cr
&
\quad-{1\over 2}\int_0^T\int_\mathbb{R}\lambda^2(t)B^2(X_t,z)\nu(dt,dz)\Big\}
+\int_0^T\int_\mathbb{R}{\mathcal  D}u(t,z)\nu(dt,dz),
\end{aligned}
$$
where ${\mathcal  D}$ is the non linear operator defined in (4.16).}

\medskip
{\it Proof.} The following well known fact will be used hereafter. If
$N_t$ ($N_0=0$) is continuous local martingale and $\big<N\big>_t$ is its predictable
quadratic variation, then the exponential process
$Z_t=\exp\big(N_t-(1/2)\big<N\big>_t
\big)$ is a continuous local martingale too, and what is more if $N'_t,N''_t$
are continuous local martingales ($N'_0=N''_0=0$) with the mutual predictable
quadratic
variation $\big<N',N''\big>_t\equiv 0$ and $Z'_t,Z''_t$ are corresponding
exponential processes, then the process $Z'_tZ''_t$ is also
local martingale which, being positive, is  a supermartingale too
(Problem 1.4.4 in \cite{23}) and so ${\mathsf E}Z'_tZ''_t\le 1, t\ge 0.$
\medskip
\noindent Let $\lambda(t)$ and $u(t,z)$ be functions involving in
Lemma 5.1. Put
$$ \begin{aligned} &
N'_t={1\over\sqrt{\varepsilon}}\int_0^t\lambda(s)B(X_s^\varepsilon,\xi_s^
\varepsilon)dW_s
\\
&
N''_t={1\over\sqrt{\varepsilon}}\int_0^tu'_z(s,\xi_s^\varepsilon)
\sigma(\xi_s^\varepsilon)dV_s.
\end{aligned}
$$
Evidently
$$
\begin{aligned}
&
\big<N'\big>_t={1\over\varepsilon}
\int_0^t\lambda^2(s)B^2(X_s^\varepsilon,\xi_s^\varepsilon)ds
\\
& \big<N''\big>_t={1\over
\varepsilon}\int_0^t(u')^2_z(s,\xi_s^\varepsilon)
\sigma^2(\xi_s^\varepsilon)ds.
\end{aligned}
\eqno(5.2)
$$
Since Wiener
processes $W_t$ and $V_t$ are independent and so
$\big<N',N''\big>_t\equiv 0$ a process $$
Z_t=\exp\Big(N'_t+N''_t-{1\over
2}\big[\big<N'\big>_t+\big<N''\big>_t \big]\Big) \eqno (5.3) $$ is
local martingale and also a supermartingale with $$ \mathsf{E}Z_t\le
1, t\ge 0. \eqno (5.4) $$ Note that $$
N'_t={1\over\varepsilon}\int_0^t\lambda(s)\big[dX_s^
\varepsilon-A(X_s^\varepsilon,\xi_s^\varepsilon)ds\big] \eqno
(5.5) $$ and also find similar representation for $N''_t$. Due to
It\^o's formula we obtain $$
{1\over\sqrt{\varepsilon}}\int_0^tu'_z(s,\xi_s^\varepsilon)
\sigma(\xi_s^\varepsilon)dV_s =u(t,\xi_t^\varepsilon)
-u(0,\xi_0)-\int_0^tu'_t(s,\xi_s^\varepsilon)ds -{1\over
\varepsilon}\int_0^t{\mathcal  L}u(s,\xi_s^\varepsilon)ds, $$ where
$\mathcal {L}=b(z){\partial d\over \partial z}+{\sigma^2(z)\over 2}
{\partial^2\over z^2}$ and consequently $$
N''_t=u(t,\xi_t^\varepsilon)-
u(0,\xi_0)+\int_0^tu'_t(s,\xi_s^\varepsilon)ds -{1\over
\varepsilon}\int_0^t{\mathcal  L}u(s,\xi_s^\varepsilon)ds. \eqno (5.6)
$$

\medskip
\noindent
(5.4) implies an obvious the inequality
$$
{\mathsf E}I(r_T(X^\varepsilon,X)+\rho_T(\nu^\varepsilon,\nu)\le\delta)
Z_T\le 1.
\eqno (5.7)
$$
which can be sharpen by changing of $Z_T$ by its lower bound. To this end
evaluate from below $\log Z_T$ on the
set $\{r_T(X^\varepsilon,X)+\rho_T(\nu^\varepsilon,\nu)\le\delta\}$.
For both
$N'_T-{1\over 2}\big<N'\big>_T$ and $N''_T-{1\over 2}\big<N''\big>_T$ we get
$$
\begin{aligned}
&
N'_T-{1\over 2}\big<N'\big>_T
\\
&
\ge {1\over\varepsilon}\sum_i\lambda(t_i)[X_{T\wedge t_{i+1}}-X_{T\wedge t_i}]
\\
&
\quad-{1\over\varepsilon}\int_0^T\int_\mathbb{R}\lambda(t)A(X_t,z)\nu(dt,dz)
-{1\over 2}\int_0^T\int_\mathbb{R}\lambda^2(t)B^2(X_t,z)\nu(dt,dz)
\\&
\quad-{1\over\varepsilon}\Big\{\sum_i|\lambda(t_i)|\Big[|X^\varepsilon_
{T\wedge t_{i+1}}-
X_{T\wedge t_{i+1}}|+|X^\varepsilon_{T\wedge t_i}-X_{T\wedge t_i}|\Big]
\\
&
\quad+\int_0^T|\lambda(t)||A(X_t^\varepsilon,\xi_t^\varepsilon)-
A(X_t,\xi_t^\varepsilon)|ds
\\
&
\quad+{1\over 2}\int_0^T\lambda^2(t)|B^2(X_t^\varepsilon,\xi_t^\varepsilon)-
B^2(X_t,\xi_t^\varepsilon)|ds
\\
&
\quad+\Big|\int_0^T\int_\mathbb{R}[\lambda(t)A(X_t,z)|+\lambda^2(t)B^2(X_t,z)]
[\nu^\varepsilon-\nu](dz,dt)\Big|\Big\}
\end{aligned}(5.8)
$$
and
$$
\begin{aligned}
&
N''_T-{1\over 2}\big<N''\big>_T
\\
&
=-{1\over\varepsilon}\int_0^T[{\mathcal  L}u(s,\xi_s^\varepsilon)+
{1\over2}u^2_z(s,\xi_s^\varepsilon)]ds+u(T,\xi_T^\varepsilon)-
u(0,\xi_0)-\int_0^Tu_t(s,\xi_s^\varepsilon)
\\
&
=-{1\over\varepsilon}\int_0^T\int_\mathbb{R}{\mathcal  D}u(s,z)\nu^\varepsilon(dz,ds)
+u(T,\xi_T^\varepsilon)-
u(0,\xi_0)-\int_0^Tu_t(s,\xi_s^\varepsilon)ds
\\
&
\ge -{1\over\varepsilon}\int_0^T\int_\mathbb{R}{\mathcal  D}u(s,z)\nu(dz,ds)
\\
&
\quad-{1\over\varepsilon}\Big\{\Big|\int_0^T\int_\mathbb{R}{\mathcal  D}u(s,z)
[\nu^\varepsilon-\nu](dz,ds)\Big|
\\
&
\quad+\varepsilon |u(T,\xi_T^\varepsilon)|+
\varepsilon|u(0,\xi_0)|+\varepsilon\int_0^T|u_t(s,\xi_s^\varepsilon)|ds
\Big\}.
\end{aligned}
\eqno(5.9)
$$
The terms in the curly brackets in the right hand sides of (5.8) and (5.9)
are random variables. Nevertheless, they can be evaluated from above on the set
$\{r_T(X^\varepsilon,X)+\rho_T(\nu^\varepsilon.\nu)\le \delta\}$ by non random
quantities. Evidently
$$
\int_0^T\int_\mathbb{R}|\lambda(t)||A(X_t^\varepsilon,\xi_t^\varepsilon)-
A(X_t,\xi_t^\varepsilon)|ds\le \text{const.}T\delta
$$
and
${1\over 2}\int_0^T\lambda^2(t)|B^2(X_t^\varepsilon,\xi_t^\varepsilon)-
B^2(X_t,\xi_t^\varepsilon)|ds\le \text{const.}\delta\int_0^T[1+|X_t|]ds.$
Denote by $H(s,z)=\lambda(s)A(X_s,z)+{\lambda^2(s)B^2(X_s,z)\over 2}\big]$.
Since $\lambda(s)$ is piece wise constant function without loss of a generality
one can assume that it is simply constant. Then function $H(s,z)$ is bounded
continuous function and so, by Lemma A.1 (see Appendix) for any $\gamma>0$
and $k\ge 1$ there exist increasing continuous
function $h^\gamma_k(y), y\ge 0$ with
$h^\gamma_k(0)=0$ and decreasing sequence $\varphi_k, k\ge1$ with
$\lim_k\varphi_k=0$ both dependent on $H(s,z)$ and $\nu$ only such that
$$
\Big|\int_0^T\int_\mathbb{R}H(s,z)[\nu^\varepsilon-\nu](ds,dz)\Big|
\le \gamma+h^\gamma_k(\delta)+\varphi_k.
$$
Further, by the remark to Lemma A.1
$$
\Big|\int_0^T\int_\mathbb{R}{\mathcal  D}u(t,z)[\nu^\varepsilon-\nu](dt,dz)\Big|
\le \gamma+h^\gamma(\delta),
$$
where $h^\gamma(y)$ is an increasing continuous function with $h^\gamma(0)=0$
depending on ${\mathcal  D}u(s,z)$ and $\nu$ only.
\medskip
\noindent
Hence, we  the lower bounds (with positive const.'s):
$$
\begin{aligned}
N'_T-{1\over 2}\big<N'\big>_T
&
\ge {1\over\varepsilon}\Big[\sum_i\lambda(t_i)
[X_{T\wedge t_{i+1}}-X_{T\wedge t_i}]-
\int_0^T\int_\mathbb{R}\lambda(t)A(X_t,z)\nu(dt,dz)
\\
&
\quad-\int_0^T\int_\mathbb{R}\lambda^2(t)B^2(X_t,z)\nu(dt,dz)
\Big]
\\
&
\quad-{\text{const.}\over\varepsilon}\Big(\gamma+h^\gamma_k(\delta)+\varphi_k
)\Big)
\end{aligned}
\eqno(5.10)
$$
and
$$
N''_T-{1\over 2}\big<N''\big>_T
\ge -{1\over\varepsilon}\int_0^T\int_\mathbb{R}{\mathcal  D}u(t,z)\nu(dt,dz)-
{\text{const.}\over\varepsilon}\Big(\varepsilon+
\gamma+h^\gamma(\delta)+\varphi_k)\Big).
\eqno (5.11)
$$
By virtue of (5.10) and (5.11) one can choose a non random lower bound:
$$
\begin{aligned}
\log Z_T
&
\ge {1\over\varepsilon}\Big[\sum_i\lambda(t_i)
[X_{T\wedge t_{i+1}}-X_{T\wedge t_i}]-
\int_0^T\int_\mathbb{R}\lambda(t)A(X_t,z)\nu(dt,dz)
\\
&
\quad-\int_0^T\int_\mathbb{R}\lambda^2(t)B^2(X_t,z)\nu(dt,dz)
-{1\over\varepsilon}\int_0^T\int_\mathbb{R}{\mathcal  D}u(t,z)\nu(dt,dz)\Big]
\\
&
\quad-{\text{const.}\over\varepsilon}\Big(\varepsilon+
\gamma+h^\gamma(\delta)+h^\gamma_k(\delta)+\varphi_k\Big).
\\
&
=\log Z_*.
\end{aligned}
$$
Hence and from (5.7), with replacing of $Z_T$ on $Z_*$, it follows
$$
\begin{aligned}
\varepsilon\log{\mathsf P}
&
\big(r_T(X^\varepsilon,X)+\rho(\nu^\varepsilon,\nu)\le
\delta\big)
\\
&
\le -\Big[\sum_i\lambda(t_i)[X_{T\wedge t_{i+1}}-X_{T\wedge t_i}]-
\int_0^T\int_\mathbb{R}\lambda(t)A(X_t,z)\nu(dt,dz)
\\
&
\quad-{1\over 2}\int_0^T\int_\mathbb{R}\lambda^2(t)B^2(X_t,z)\nu(dt,dz)\Big]
+\int_0^T\int_\mathbb{R}{\mathcal  D}u(t,z)\nu(dt,dz)
\\
&
\quad+\text{const.}\Big\{\big(\varepsilon+\gamma+h^\gamma(\delta)
+h^\gamma_k(\delta)+\varphi_k
)\Big\},
\end{aligned}
\eqno(5.12)
$$

\medskip
\noindent
The desired result holds since the term in the curly brackets of the right
hand side in (5.12) goes to zero if limit
`$lim_{k\to\infty}lim_{\gamma\to 0}lim_{\delta\to 0}\lim_
{\varepsilon\to 0}$' is taken.

\medskip
\noindent
{\it Proof of Theorem 5.1.} follows from Lemmas 5.1, A.2, and A.3
(see Appendix) since
$$
\begin{aligned}
&
-\sup_\lambda\Big\{\sum_i\lambda(t_i)[X_{T\wedge t_{i+1}}-X_{T\wedge t_i}]-
\int_0^T\int_\mathbb{R}\lambda(t)A(X_t,z)\nu(dt,dz)
\\
&
-{1\over 2}\int_0^T\int_\mathbb{R}\lambda^2(t)B^2(X_t,z)\nu(dt,dz)\Big\}
+\inf_u\int_0^T\int_\mathbb{R}{\mathcal  D}u(t,z)\nu(dt,dz)
\\
&
=-\Big[{1\over 2}S_T(X,\nu)+{1\over 8}F_T(\nu)\Big]=-L_T(X,\nu).
\end{aligned}
$$

\medskip
\noindent\vskip .25in
\noindent{\bf 6. Lower bound for local LDP in $\Bbb{C}_{[0,T]}\times
\Bbb{M}_{[0,T]}$}

\medskip
{\bf Theorem 6.1.} {\it Under {\rm (A.1), (A.2)}, and {\rm (A.3)}, for every
$(X,\nu)$ from $\Bbb{C}_{[0,T]}\times \in\Bbb{M}_{[0,T]}$
$$
\underline{\lim}_{\delta\to 0}\underline{\lim}_{\varepsilon\to 0}
\varepsilon\log {\mathsf P}\big(r_T(
X^\varepsilon,X)+\rho_T(\nu^\varepsilon,\nu)\le \delta\big)\ge -L_T(X,\nu).
$$}

\medskip
\noindent
Evidently for $X,\nu$ such that $L_T(X,\nu)=\infty$ it is nothing
to prove. Therefore below we consider only the case
$L_T(X,\nu)<\infty$ which distinguishes subsets from
$\Bbb{C}_{[0,T]}\times\Bbb{M}_{[0,T]}$:

{\bf (i)} $dX_t\ll dt$ and ${1\over 2}S_T(X,\nu)=
{1\over 2}\int_0^T{[\dot{X}_t-A_\nu(t,X_t)]^2\over B_\nu^2(t,X_t)}dt<\infty$;

{\bf (ii)} $d\nu=nd\lambda, \ d_zn=n'_zdz$ and
${1\over 8}F_T(\nu)={1\over 2}\int_0^T\int_\mathbb{R}{v_\nu^2(t,z)
\over \sigma^2(z)}n(t,z)dzdt<\infty$,
where
$$
v_\nu(t,z)={\sigma^2(z)\over 2}\Big[{n'_z(t,z)\over n(t,z)}-{p'(z)\over p(z)}
\Big].
\eqno (6.1)
$$
It is convenient to consider further subset (ii') of (ii):

{\bf (ii')} the function $v_\nu(t,z)$
is compactly supported in $z$ and continuously differentiable in $(t,z),$
having bounded partial derivatives.

\medskip
\noindent
The central role in proving Theorem 6.1 plays

\smallskip
{\bf Lemma 6.1.} {\it Assume {\rm (i), (ii')}, and
$\inf_{x,z}B^2(x,z)>0$. Then for any $\delta>0$ and
$\gamma>0$ there exists an increasing continuous function $h_\gamma(y)$ with
$h_\gamma(0)=0$, depending on ${v_\nu^2(s,z)\over \sigma^2(z)}$ and $\nu$
only, such that
$$
\varliminf_{\varepsilon\to 0}
\varepsilon\log {\mathsf P}\big(r_T(
X^\varepsilon,X)+\rho_T(\nu^\varepsilon,\nu)\le \delta\big)\ge -L_T(X,\nu)
-\gamma-h_\gamma(\delta).
$$}

\medskip
{\it Proof.} Put
$$
\begin{aligned}
&
b_\nu(t,z)=b(z)+v_\nu(t,z)
\\
&
G_\nu(t,x,z)={\dot{X}_t-A_\nu(t,X_t)\over B_\nu(t,X_t)}B(x,z)+A(x,z)
\end{aligned}
\eqno(6.2)
$$
and parallel to $(X_t^\varepsilon,\xi_t^\varepsilon)$ introduce, on the same
stochastic basis, new diffusion pair
$(\widetilde{X}_t^\varepsilon,\widetilde{\xi}_t^\varepsilon)$:
$$
\begin{aligned}
&
d\widetilde{X}_t^\varepsilon=G_\nu(t,\widetilde{X}_t^\varepsilon,
\widetilde{\xi}_t^\varepsilon)dt+
\sqrt{\varepsilon}B(\widetilde{X}_t^\varepsilon,
\widetilde{\xi}_t^\varepsilon)dW_t
\\
&
d\widetilde{\xi}_t^\varepsilon={1\over \varepsilon}b_\nu(t,
\widetilde{\xi}_t^\varepsilon)dt+
{1\over\sqrt{\varepsilon}}\sigma(\widetilde{\xi}_t^\varepsilon)dV_t
\end{aligned}
\eqno(6.3)
$$
subject to the same initial point $(x_0$, $\xi_0)$. Also denote by
$\widetilde{\nu}^\varepsilon(dt,dz)$ the occupation measure corresponding to
$\widetilde{\xi}^\varepsilon$:
$\widetilde{\nu}^\varepsilon(\Delta\times\Gamma)=
\int_0^\infty I(t\in \Delta, \widetilde{\xi}_t^\varepsilon\in\Gamma)dt.$
\medskip
\noindent
By virtue of the formula $b(z)={1\over 2}{p'(z)\over p(z)}+
\sigma'(z)\sigma(z)$ (see (2.1)) we get ${2b_\nu(t,z)\over\sigma^2(z)}=
{n'_z(t,z)\over n(t,z)}+2{\sigma'(z)\over \sigma(z)}$ and so
$p_\nu(t,z)=c(t){\exp\big(2\int_0^z{b_\nu(t,y)\over\sigma^2(y)}dy\big)\over
\sigma^2(z)}$, with norming constant $c(t)$ such that $\int_\mathbb{R}p_\nu(t,z)=1$,
coincides with $n(t,z)$. Then by Lemma A.5 (see Appendix)
$$
\mathsf{P}\text{-}\lim_{\varepsilon\to 0}\rho_T(\widetilde{\nu}^\varepsilon,\nu)=0
\quad\text{and}\quad
\mathsf{P}\text{-}\lim_{\varepsilon\to 0}r_T(\widetilde{X}^\varepsilon,X)=0.
\eqno (6.4)
$$
\medskip
\noindent
Denote by $Q^\varepsilon$ and $\widetilde{Q}^\varepsilon$
distributions of $(X_t^\varepsilon,\xi_t^\varepsilon)_{t\le T}$,
$(\widetilde{X}_t^\varepsilon,\widetilde{\xi}_t^\varepsilon)_{t\le T}$
respectively. By Theorem 7.18 (Ch. 7 in \cite{24}
$Q^\varepsilon$ is absolutely continuous w.r.t.
$\widetilde{Q}^\varepsilon$ and
$$
{dQ^\varepsilon \over d\widetilde{Q}^\varepsilon}(\widetilde{X}^\varepsilon,
\widetilde{\xi}^\varepsilon)=\exp\Big({1\over\sqrt{\varepsilon}}M_T^
\varepsilon-{1\over 2\varepsilon}\big<M^\varepsilon\big>_T+
{1\over \sqrt{\varepsilon}}M_T-{1\over 2\varepsilon}\big<M\big>_T\Big),
\eqno (6.5)
$$
where
$$
\begin{aligned}
&
M_t^\varepsilon=-\int_0^t{v_\nu(s,\widetilde{\xi}_s^\varepsilon)
\over\sigma(\widetilde{\xi}_s^\varepsilon)}dV_s
\quad\text{and}\quad M_t=-\int_0^t
{\dot{X}_s-A_\nu(s,X_s)\over B_\nu(s,X_s)}dW_s
\\
&
\big<M^\varepsilon\big>_t=\int_0^t{v_\nu^2(s,\widetilde{\xi}_s
^\varepsilon)\over\sigma^2(\widetilde{\xi}_s^\varepsilon)}ds
\quad\text{and}\quad
\big<M\big>_t=\int_0^t{[\dot{X}_s-A_\nu(s,X_s)]^2\over B^2_\nu(s,X_s)}ds.
\end{aligned}
$$
By virtue of (6.5) and the rule of changing for probability measure we
obtain
$$
{\mathsf P}\big(r_T(X^\varepsilon,X)+\rho_T(\nu^\varepsilon,\nu)\le \delta\big)=
{\mathsf E}\Big[{dQ^\varepsilon \over d\widetilde{Q}^\varepsilon}
(\widetilde{X}^\varepsilon,\widetilde{\xi}^\varepsilon)
I(r_T(\widetilde{X}^\varepsilon,X)+\rho_T(\widetilde{\nu}
^\varepsilon,\nu)\le \delta)\Big].
\eqno (6.6)
$$
The desired lover bound can be derived from (6.6) provided that a relevant
lover bound for the right hand side of (6.6) can be found. Use an  obvious
inequality:
$$
I(r_T(\widetilde{X}^\varepsilon,X)+\rho_T(\widetilde{\nu}
^\varepsilon,\nu)\le \delta)\ge
I(r_T(\widetilde{X}^\varepsilon,X)+\rho_T(\widetilde{\nu}
^\varepsilon,\nu)\le \delta, \ |M_T^\varepsilon|\le k, \ |M_T|\le k)
$$
and estimate from below $\log {dQ^\varepsilon \over d\widetilde{Q}
^\varepsilon}(\widetilde{X}^\varepsilon,\widetilde{\xi}^\varepsilon)$
on the set $\{r_T(\widetilde{X}^\varepsilon,X)+
\rho_T(\widetilde{\nu}
^\varepsilon,\nu)\le \delta\}\cap\{|M_T^\varepsilon|\le k\}\cap
\{|M_T|\le k\}$.
Noticing that ${1\over 2}\big<M\big>_T={1\over 2}S(X,\nu)_T$ and
$$
{1\over 2}\int_0^T\int_\mathbb{R}{v_\nu^2(s,z)\over\sigma^2(z)}n(t,z)dzdt=
{1\over 8}F_T(\nu)
$$
 we obtain
$$
\log {dQ^\varepsilon \over d\widetilde{Q}
^\varepsilon}(\widetilde{X}^\varepsilon,\widetilde{\xi}^\varepsilon)
\ge -{2k\over \sqrt{\varepsilon}}-{1\over\varepsilon}L_T(X,\nu)-
{1\over 2\varepsilon}\Big|\int_0^T\int_\mathbb{R}
{v_\nu^2(s,z)\over\sigma^2(z)}[\widetilde{\nu}^\varepsilon(dt,dz)
-n(t,z)dzdt\Big|.
$$
By the remark to Lemma A.1 (see Appendix), for any $\gamma>0$ there
exists increasing continuous
function $h_\gamma(y)$ with $h\gamma(0)=0$, depending on
${v_\nu^2(s,z)\over\sigma^2(z)}$ and $\nu$ only, such that
$$
{1\over 2}\Big|\int_0^T\int_\mathbb{R}
{v_\nu^2(s,z)\over\sigma^2(z)}[\widetilde{\nu}^\varepsilon(dt,dz)
-n(t,z)dzdt\Big|\le \gamma+h_\gamma(\delta).
$$
Then the lower bound for the right hand side of (6.6) is the
following: $-{2k\over \sqrt{\varepsilon}}-
{1\over\varepsilon}[L_T(X,\nu)+\gamma+h_\gamma(\delta)]$. It implies
$$
\begin{aligned}
\varepsilon\log {\mathsf P}\big(r_T(
X^\varepsilon,X)
&
+\rho_T(\nu^\varepsilon,\nu)\le \delta\big)
\ge -L_T(X,\nu)-2k\sqrt{\varepsilon}-\gamma-h_\gamma(\delta)
\\
&
\quad+\varepsilon\log {\mathsf P}\big(r_T(\widetilde{X}
^\varepsilon,X)+\rho_T(\widetilde{\nu}^\varepsilon,\nu)\le \delta,
\ |M_T^\varepsilon|\le k, \ |M_T|\le k\big).
\end{aligned}
$$
\medskip
\noindent
Thus the statement of the lemma holds since
$$
\lim_k\lim_{\varepsilon\to 0}{\mathsf P}\big(r_T(\widetilde{X}
^\varepsilon,X)+\rho_T(\widetilde{\nu}^\varepsilon,\nu)\le \delta,
\ |M_T^\varepsilon|\le k, \ |M_T|\le k\big)=1
$$
what follows by virtue of (6.4), obvious $\lim_k\mathsf{P}(|M_T|>k)=0$
and
$$
\mathsf{P}(|M_T^\varepsilon|>k)\le {\mathsf{E}|M_T^\varepsilon|^2\over k^2}
={\mathsf{E}\big<M^\varepsilon\big>_T\over k^2}\le{\text{const.}\over k^2}\to 0,
k\to\infty.
$$

\medskip
\noindent
{\it Proof of Theorem 6.1.} Assume (i), (ii), and $\inf_{x,z}B^2(x,z)>0$.
Due to Lemma A.4 (see Appendix), one can choose a sequence $\nu^{(k)}, k\ge 1$ of measures
such that for every  $k$ the function $v_{\nu^{(k)}}(t,z)$ satisfies
(ii') and what is more $\rho(\nu,\nu^{(k)})\to 0$, $L_T(X,\nu^{(k)})\to
L_T(X,\nu)$.  On the other hand, by Lemma 6.1
for any $\delta>0$ and
$\gamma>0$ there exist increasing continuous function $h_{\gamma,k}(y)$ with
$h_{\gamma,k}(0)=0$, depending on $\nu^{(k)}$, such that
$$
\underline{\lim}_{\varepsilon\to 0}
\varepsilon\log {\mathsf P}\big(r_T(
X^\varepsilon,X)+\rho_T(\nu^\varepsilon,\nu^{(k)})\le \delta\big)
\ge -L_T(X,\nu^{(k)})-\gamma-h_{\gamma,k}(\delta).
$$
Choose $k_\circ(\delta)$ such that for any $k\ge k_\circ(\delta)$ we have
$0<\delta-\rho_T(\nu,\nu^{(k)})\le \delta/2$. Then, taking into account
the triangular inequality:
$\rho_T(\nu^\varepsilon,\nu)\le \rho_T(\nu^\varepsilon,
\nu^{(k)})+\rho_T(\nu,\nu^{(k)})$, we arrive at a lower bound:
$$
\begin{aligned}
\underline{\lim}_{\varepsilon\to 0}
\varepsilon\log {\mathsf P}\big(r_T(
X^\varepsilon,X)
&+\rho_T(\nu^\varepsilon,\nu)\le \delta\big)
\\
&
\ge \underline{\lim}_{\varepsilon\to 0}
\varepsilon\log {\mathsf P}\big(r_T(
X^\varepsilon,X)+\rho_T(\nu^\varepsilon,\nu^{(k)})\le \delta/2\big)
\\
&
\ge -L_T(X,\nu^{(k)})-\gamma-h_{\gamma,k}(\delta/2).
\end{aligned}
$$
The right hand side of the last inequality converges
to $-L_T(X,\nu)$ if limit
\newline
`$\lim_k \lim_{\gamma\to 0} \lim_{\delta\to 0}$' is taken.
\medskip
\noindent
Assume only (i) and (ii). Parallel to the process $X_t^\varepsilon$
introduce new diffusion $X_t^{\varepsilon,\beta},\beta\neq 0$:
$$
dX_t^{\varepsilon,\beta}=A(X_t^{\varepsilon,\beta},\xi_t^\varepsilon)dt+
\sqrt{\varepsilon}\big[B(X_t^{\varepsilon,\beta},\xi_t^\varepsilon)dW_t+
\beta dW'_t\big]
$$
subject to the same initial point $x_0$, where $W'_t$ is a Wiener process
independent of $(W_t, \xi_t^\varepsilon)$. The diffusion parameter here
is $B^2(x,z)+\beta^2$ and so, due to proved above,
$$
\underline{\lim}_{\delta\to 0}\underline{\lim}_{\varepsilon\to 0}
\varepsilon\log {\mathsf P}\big(r_T(
X^{\varepsilon,\beta},X)+\rho_T(\nu^\varepsilon,\nu)
\le \delta\big)\ge -L_T^\beta(X,\nu),
$$
where $L_T^\beta(X,\nu)={1\over 2}S^\beta_T(X,\nu)+{1\over 8}F_T(\nu)$,
and where
$$
S^\beta_T(X,\nu)=\int_0^\infty{[\dot{X}_t-A_\nu(y,X_t)]^2
\over B^2_\nu(t,X_t)+\beta^2}dt.
$$
Evidently $\lim_{\beta\to 0}S_T^\beta(X,\nu)=S_T(X,\nu)$.
On the other hand, by Lemma A.6 (see Appendix)
$$
\lim_{\beta\to 0}\bar{\lim}_{\varepsilon\to 0}\varepsilon\log
\mathsf{P}(r_T(X^{\varepsilon,\beta},X^\varepsilon)>\eta)=-\infty.
$$
To get the desired result, we combine both these facts. Namely, using
the triangular
inequality: $r_T(X^\varepsilon,X)\le
r_T(X^{\varepsilon,\beta},X^\varepsilon)+
r_T(X^{\varepsilon,\beta},X)$ and taking $\eta=\delta/2$ we arrive at
an upper bound
$$
\begin{aligned}
&
{\mathsf P}\big(r_T(
X^{\varepsilon,\beta},X)+\rho_T(\nu^\varepsilon,\nu)\le \delta\big)
\\
&
\quad \le {\mathsf P}\big(r_T(
X^\varepsilon,X)+\rho_T(\nu^\varepsilon,\nu)\le \delta/2\big)
\\
&
\quad +{\mathsf P}\big(r_T(X^{\varepsilon,\beta},X^\varepsilon)>\delta/2\big)
\\
&
\quad
\le2\max\Big[{\mathsf P}\big(r_T(
X^\varepsilon,X)+\rho_T(\nu^\varepsilon,\nu)\le \delta/2\big),
{\mathsf P}\big(r_T(X^{\varepsilon,\beta},X^\varepsilon)>\delta/2\big)\Big]
\end{aligned}
$$
which implies
$$
\underline{\lim}_{\delta\to 0}\underline{\lim}_{\varepsilon\to 0}
\varepsilon\log {\mathsf P}\big(r_T(
X^\varepsilon,X)+\rho_T(\nu^\varepsilon,\nu)\le \delta\big)\ge -
\lim_{\beta\to 0}L^\beta_T(X,\nu)=-L_T(X,\nu).
$$
\medskip
\noindent
Other approach for establishing lower bound with singular diffusion parameter
can be found in Puhalskii \cite{25}.

\section{\bf Proof of main result}

{\it Proof of Theorem 3.1.} Due to Theorems 4.1, 5.1, and Proposition 3.1
the family $(X^\varepsilon,\nu^\varepsilon)$ obeys the LDP in
$(\Bbb{C}_{[0,k]}\times\Bbb{M}_{[0,k]}, r_k\times\rho_k)$ with rate
function $L_k(X,\nu)$. Then it obeys the LDP in the metric space
$(\Bbb{C}\times\Bbb{M}, r\times\rho)$ with rate function
$\sup_kL_k(X,\nu)=L(X,\nu)$.
\medskip
\noindent
{\it Proof of Corollary 2.1.} The result holds since
$\inf_{X\in\Bbb{C}}S(X,\nu)$ is attained at $X_t^\circ$, being a solution of
a differential equation: $\dot{X}_t=A_\nu(t,X_t^\circ)$ subject to
$X_0^\circ=x_0$, and so $S(X^\circ,\nu)=0$.

\medskip
\noindent
{\it Proof of Corollary 2.2.} The first statement is obvious.
\medskip
\noindent
Assume $B^2(x,z)\equiv 0$. In this case $S(X,\nu)=0$ for any
$X_t$ being a solution of a differential equation
$\dot{X}_t=\int_\mathbb{R}A(X_t,z)n(t,z)dz$ subject to $X_0=x_0$; otherwise
$S(X,\nu)=\infty$. Therefore
$$
L(X,\nu)=\begin{cases}
{1\over 8}\inf_{\nu:\dot{X}_t=\int_\mathbb{R}A(X_t,z)n(t,z)dz, X_0=x_0}F(\nu)\\
\infty, & \text{otherwise}.
\end{cases}
$$
On the other hand, since $F(\nu)<\infty$ implies
$d\nu=nd\lambda, \ d_zn=n'_zdz$,
assuming measurability in $t$ of function
$$
H(t,\dot{X}_t,X_t)=\inf_{\nu:\dot{X}_t=\int_\mathbb{R}A(X_t,z)n(t,z)dz, X_0=x_0}
\int_\mathbb{R}\sigma^2(z)\Big[
{n'_z(t,z)\over n(t,z)}-{p'(z)\over p(z)}\Big]^2n(t,z)dz
\eqno (7.1)
$$
we arrive at independent of $t$ function $H(t,y,x)\equiv H(y,x)$, or by other
words, `$inf$' in (7.1) can be taken over all measures $\nu$ with densities
$n(t,z)\equiv m(z)$. The last means the desired result holds if the function
$$
H(y,x)=\inf_{m:\begin{cases}
dm=m'dz
\\
y=\int_\mathbb{R}A(x,z)m(z)dz
\end{cases}}
\int_\mathbb{R}\sigma^2(z)\Big[
{m'(z)\over m(z)}-{p'(z)\over p(z)}\Big]^2m(z)dz
\eqno (7.2)
$$
is measurable. We check this by showing that level sets of $H(y,x)$ are closed.
\medskip
\noindent
Let $c\ge 0$ be fixed and $(y_n,x_n),n\ge 1$ be a sequence from
$\{(y,x):H(y,x)\le c\}$ converging to a limit point
$(y_0,x_0)$. Show that $H(y_0,x_0)\le c$. By virtue of assumption
(A.1) the set ${\mathcal  A}(y,x)=\{m:y=\int_\mathbb{R}A(x,z)m(z)dz\}$ is closed in the
Levy-Prohorov metric that is for every fixed $(y,x)$ there
exists a density $m^{(y,x)}$ from ${\mathcal  A}(y,x)$ such that
$$
H(y,x)=\begin{cases}
\int_\mathbb{R}\sigma^2(z)\Big[
{(m^{(y,x)}(z))'\over m^{(y,x)}(z)}-{p'(z)\over p(z)}\Big]^2m^{(y,x)}(z)dz,
& dm^{(y,x)}=(m^{(y,x)})'dz
\\
\infty, &\text{otherwise}.
\end{cases}
\eqno (7.3)
$$
Note that the function $H(y,x)$, defined in (7.3), obeys a following property:
there exists a measure $\nu^{(y,x)}$ from $\Bbb{M}_{[0,1]}$, having
density $m^{(y,x)}(z)$ w.r.t. $dtdz$, such that $H(y,x)=F_1(\nu^{(y,x)})$.
Since ${1\over 8}F_1(\nu)$ is good rate function level sets
$\{y,x:H(y,x)\le c\}$ are compacts. Therefore $H(y_0,x_0)\le c.$

\section{\bf Appendix}

\noindent{\bf 1.} $\underline{\text{Evaluation via Levy-Prohorov's metric.}}$

\smallskip
{\bf Lemma A.1.} {\it Let $T>0$, $\nu',\nu''\in \Bbb{M}_{[0,T]}$,
$\rho_T(\nu',\nu'')=q$,
and $f(t,z)$ be bounded continuous function.
Then
for any $\gamma>0$ and $k\ge 1$ one can choose increasing continuous function
$h^\gamma_k(y), y\ge 0$ with $h^\gamma_k(0)=0$ and decreasing sequence
$\varphi_k, k\ge 1$ with $\lim_k\varphi_k=0$ both depending on $f(t,z)$
and only from one of $\nu'$ or $\nu''$ such that
$$
\Big|\int_0^T\int_\mathbb{R}f(t,z)[\nu'-\nu''](dt,dz)\Big|\le \gamma+h^\gamma_k(q)
+\varphi_k.
$$}

\smallskip
{\bf Remark} {\it If $f(t,z)$ is bounded compactly
supported continuous function, then the statement of the lemma
remains true with $h^\gamma_k(y)\equiv h^\gamma(y)$ and
$\varphi_k\equiv 0.$}

\medskip
{\it Proof.} Assume $f(t,z)$ is continuously
differentiable (one in $z$ and twice in $(t,z)$) and compactly
supported in $z$. Denote by $F'(t,z)=\nu'([0,t]\times(-\infty,z])$
that is $F'(t,z)$ is the distribution function corresponding to
$\nu'$. Integrating by parts we get
$$
\int_0^T\int_\mathbb{R}f(t,z)\nu'(dt,dz)=-
\int_0^T\int_\mathbb{R}\Big[{\partial f(t,z)\over\partial z}+{\partial^2 f(t,z)
\over\partial t\partial z}\Big]F'(t,z)dzdt
$$
and consequently ($F''$ is the  distribution functions corresponding to
$\nu''$)
$$
\Big|\int_0^T\int_\mathbb{R}f(t,z)[\nu'-\nu''](dt,dz)\Big|\le
\int_0^T\int_\mathbb{R}|F'(t,z)-F''(t,z)|m(t,z)dzdt,
$$
where $m(t,z)=|{\partial\over \partial z}f(t,z)|+|{\partial^2\over
\partial t\partial z}f(t,z)|$.

\medskip
\noindent
Assume $f(t,z)$ is compactly supported in $z$ and continuous only. Then,
approximating it by  compactly  supported  and continuously
differentiable in $z$ function $f^\gamma(t,z)$ in a sense
$
\sup_{t,z}|f(t,z)-f^\gamma(t,z)|\le {\gamma\over 2T},
$
due the foregoing proof, we get
$$
\Big|\int_0^T\int_\mathbb{R}f(t,z)[\nu'-\nu''](dt,dz)\Big|\le \gamma+
\int_0^T\int_\mathbb{R}|F'(t,z)-F''(t,z)|m^\gamma(t,z)dzdt
$$
with $m^\gamma(t,z)=|{\partial\over \partial z}f^\gamma(t,z)|+
|{\partial^2\over\partial t\partial z}f^\gamma(t,z)|$

\medskip
\noindent
In the general case, one can choose a decomposition $f(t,z)=f_k(t,z)+
g_k(t,z)$, where $f_k(t,z)$ is continuous compactly supported in $z$ on
the interval
$[-k,k]$ function while $g_k(t,z)\equiv 0$ on the interval $[-(k-1/2),(k-1/2)]$
and is bounded: $|g_k(t,z)|\le L.$
Then by foregoing result we get
$$
\begin{aligned}
\Big|\int_0^T\int_\mathbb{R}f(t,z)[\nu'-\nu''](dt,dz)\Big|\le \gamma
&
+\int_0^T\int_\mathbb{R}|F'(t,z)-F''(t,z)|m^\gamma_k(t,z)dzdt
\\
&+L\int_0^T\int_{|z|>k-1/2}[\nu'+\nu''](dt,dz),
\end{aligned}
$$
where $m^\gamma_k(t,z)=|{\partial\over \partial z}f^\gamma_k(t,z)|+
|{\partial^2\over\partial t\partial z}f^\gamma_k(t,z)|$. Evaluate
from above the last integral from the right hand side.
To this end, choose an increasing sequences $z_k\nearrow\infty$,
$k\to\infty$ such that $z_k\le k-1/2$ and for every $k$  $z_k$
and $-z_k$ are points of continuity for the distribution function $F'(T,z)$.
Then
$$
\begin{aligned}
\int_0^T\int_{|z|>k-1/2}[\nu'+\nu''](dt,dz)
&
\le 2\int_0^T\int_{|z|>z_k}\nu'(dt,dz)
\\
&
\quad +\Big|\int_0^T\int_{|z|>z_k}[\nu'-\nu''](dt,dz)\Big|
\\
&
\le 2\int_0^T\int_{|z|>z_k}\nu'(dt,dz)
\\
&
\quad +\big|F'(T,z_k)-F''(T,z_k)\big|
\\
&
\quad +\big|F'(T,-z_k)-F''(T,-z_k)\big|.
\end{aligned}
$$

\medskip
\noindent
Now evaluate from above
$|F'(t,z)-F''(t,z)|$ via $q$ and $F'(t,z)$. From the
definition of the Levy-Prohorov metric (see e.g. \cite{22}, \cite{26})
it follows:
$q+F'(t-q,z-q)-F'(t,z)\le F'(t,z)-F''(t,z)\le
q+F'(t+q,z+q)-F'(t,z)$
and so
$$
|F'(t,z)-F''(t,z)|\le q+[F'(t+q,z+q)-F'(t-q,z-q)].
$$

\medskip
\noindent
Hence, combining all obtained upper estimates, we arrive at the desired
result with
$$
\begin{aligned}
h^\gamma_k(y)
&
=y\Big(2L+\int_0^T\int_\mathbb{R}m^\gamma_k(t,z)dtdz\Big)
\\
&
\quad+\int_0^T\int_\mathbb{R}[F'(t+y,z+y)-F'(t-y,z-y)]m^\gamma_k(t,z)dzdt+
\\
&
\quad +L\big|F'(T+y,z_k+y)-F'(T-y,z_k-y)\big|
\\
&
\quad +L\big|F'(T+y,-z_k+y)-F'(T-y,-z_k-y)\big|
\end{aligned}
$$
and
$$
\varphi_k=2L\int_0^T\int_{|z|>z_k}\nu'(dt,dz).
$$
\medskip
\noindent
The same proof takes place with $F''$ instead of $F'$.

\medskip
\noindent
{\bf 2.} $\underline{\text{The Fenchel-Legendre transform.}}$
\newline Let $\lambda(t)=\sum_i\lambda(t_i)I(t_i\le t
<t_{i+1})$ with non overlapping intervals $[t_i,t_{i+1})$.
For any  $X\in \Bbb{C}_{[0,T]}$ and $\nu\in\Bbb{M}_{[0,T]}$
put $\int_0^T\lambda(t)dX_t=
\sum_i\lambda(t_i)[X_{T\wedge t_{i+1}}-X_{T\wedge t_i}]$,
$A_\nu(t,X_t)=\int_\mathbb{R}A(X_t,z)K_\nu(t,dz)$, and
$B^2_\nu(t,X_t)=\int_\mathbb{R}B^2(X_t,z)K_\nu(t,dz)$.
Let $\mathcal {D}$ be non linear operator defined in (4.16).

\smallskip
{\bf Lemma A.2.} {\it For any  $X\in \Bbb{C}_{[0,T]}$ and
$\nu\in\Bbb{M}_{[0,T]}$
$$
\begin{aligned}
\sup\int_0^T
&
\big[\lambda(t)dX_t-(A_\nu(t,X_t)-{1\over 2}
\lambda^2(t)B^2_\nu(t,X_t))dt
\\
&=
\begin{cases}
{1\over 2}\int_0^T {[\dot{X}_t-A_\nu(t,X_t)]^2\over B^2_\nu(t,X_t)}dt &
dX_t=\dot{X}_tdt
\\
\infty, & \text{otherwise},
\end{cases}
\end{aligned}
$$
where `sup' is taken over all piece wise
constant functions $\lambda(t)$.}

\smallskip
{\bf Lemma A.3.} {\it For any $\nu\in\Bbb{M}_{[0,T]}$
$$
\begin{aligned}
\inf\int_0^T
&
\int_\mathbb{R}\mathcal {D}u(t,z)\nu(dt,dz)
\\
&
=
\begin{cases}
-{1\over 8}\int_0^T\int_\mathbb{R}\sigma^2(z)\Big[{n'_z(t,z)\over n(t,z)}-
{p'(z)\over p(z)}\Big]^2n(t,z)dzdt, & d\nu=nd\lambda, \ d_zn=n'_zdz
\\
-\infty, &\text{otherwise},
\end{cases}
\end{aligned}
$$
where `$inf$' is taken over all
twice continuously differentiable (once in $t$ and twice in $z$)
compactly supported in $z$ functions $u(t,z).$}

\medskip
{\it Proof of Lemma A.2.} For $dX_t\not\ll dt$ the result follows from
Lemma 6.1 in \cite{17} (see also Lemma 2.1 in \cite{27}). For
$dX_t=\dot{X}_tdt$ by lemma 6.1 \cite{17} `$sup\int_0^T$'
is equal $\int_0^T\sup_{\lambda\in \mathbb{R}}\Big\{\lambda\Big(\dot{X}_t-
A_\nu(t,X_t)\Big)-{1\over 2}\lambda^2B^2_\nu(t,X_t)\Big)\Big\}dt
={1\over 2}\int_0^T{[\dot{X}_t-A_\nu(t,X_t)]\over B^2_\nu(t,X_t)}dt$.

\medskip
\noindent
{\it Proof of Lemma A.3.} Assume $d\nu=nd\lambda, \ d_zn=n'_zdz$.
Due to (2.1), ${p'(z)\over p(z)}=
{2b(z)-2\sigma(z)\sigma'(z)\over \sigma^2(z)}$ and so
$b(z)={1\over 2}\big[\sigma^2(z){p'(z)\over p(z)}+2\sigma(z)\sigma'(z)\big].$
Putting $v(t,z)=u'_z(t,z)$ and taking into account the formula for $b(z)$
we arrive at
$$
\begin{aligned}
\int_0^T\int_\mathbb{R}{\mathcal  D}u(t,z)n(t,z)dzdt
&
={1\over 2}\int_0^T\int_\mathbb{R}\Big\{\big[\sigma^2(z){p'(z)\over p(z)}
+2\sigma(z)\sigma'(z)\big]v(t,z)
\\
&
\quad +\sigma^2(z)(v'_z(t,z)+v^2(t,z))\Big\}n(t,z)dzdt.
\end{aligned}
\eqno(8.1)
$$
Then, integrating by parts:
$$
\int_\mathbb{R}\sigma^2(z)(v'_z(t,z)n(t,z)dz=
-\int_\mathbb{R}v(t,z)[2\sigma(z)\sigma'(z)n(t,z)+\sigma^2(z)n'_z(t,z)]dz,
$$
we obtain
$$
\begin{aligned}
&
\int_0^T\int_\mathbb{R}{\mathcal  D}u(t,z)n(t,z)dzdt
\\
&
={1\over 2}\int_0^T\int_\mathbb{R}\sigma^2(z)\Big(v^2(t,z)n(t,z)
+v(t,z)\big[{p'(z)\over p(z)}
n(t,z)-n_z(t,z)\big]\Big)dzdt.
\end{aligned}
\eqno(8.2)
$$
(8.2) and the method of proving for lemma 6.1 in \cite{17} imply
$$
\begin{aligned}
\inf\int_0^T\int_\mathbb{R}{\mathcal  D}u(t,z)n(t,z)dzdt
&
={1\over 2}\int_0^T\int_\mathbb{R}\sigma^2(z)\inf_{v\in \mathbb{R}}\Big(v^2n(t,z)
\\
&\quad+v\big[{p'(z)\over p(z)}n(t,z)-n_z(t,z)\big]\Big)dzdt
\\
&
=-{1\over 8}\int_0^T\int_\mathbb{R}\sigma^2(z)
\big[{n_z(t,z)\over n(t,z)}-{p'(z)\over p(z)}\big]^2n(t,z)dzdt.
\end{aligned}
$$
Thus for `$d\nu=nd\lambda, d_zn=n'_zdz$', the result holds.

\medskip
\noindent
Assume $d\nu=n\lambda, d_zn\not\ll dz$. Show that
$\inf\int_0^T\int_\mathbb{R}{\mathcal  D}u(t,z)n(t,z)dzdt=-\infty$.  To this end, take
$u(t,z)\equiv u(z)$
and put $v(z)=u'(z)$. The function $v(z)$ is compactly supported and
continuously differentiable and, in particular, has the finite total variation.
Put $n(z)=\int_0^Tn(t,z)dt$ and $w(z)={1\over 2}\sigma^2(z)n(z)$.
It is clear that there exists a positive constant, say, $\ell$ such that
$I(v)=\ell\int_\mathbb{R}[v^2(z)+|v(z)|]n(z)dz+\int_\mathbb{R}w(z)dv(z)$ is an upper bound
for  the right hand side of (8.1). Show that $I(v)$ can be chosen less
than any
negative quantity. Use the fact that $I(v)$ is well defined
not only for compactly supported and continuously differentiable function
$v(z)$ but also for any compactly supported function $v^\alpha(z)$ obeying
finite total variation. Assume that there exists a family of $v^\alpha(z),
\alpha\in (0,1]$ such that
$$
\lim_{\alpha\to 0}I(v^\alpha)=-\infty
\eqno (8.3)
$$
and every function $v^\alpha(z)$ obeys an approximation
by $v^\alpha_m(z), m\ge 1$ of continuously differentiable compactly supported
functions in a sense
$$
\lim_mI(v^\alpha_m)=I(v^\alpha).
\eqno (8.4)
$$
We show that under (8.3) and (8.4) the desired result holds. In fact,
for fixed $\alpha$ one can choose a number $m_\alpha$ such that
$|I(v^\alpha)-I(v^\alpha_{m_\alpha})|\le 1$. Hence we obtain
$$
\inf\int_0^T\int_\mathbb{R}{\mathcal  D}u(t,z)n(t,z)dzdt\le I(v^\alpha_{m_\alpha})\le
1+I(v^\alpha)\to -\infty, \ \alpha\to 0.
$$
Therefore, only (8.3) and (8.4) have to be checked.
Since $d_zn\not\ll dz$ the function $n(z)$
is not absolutely continuous and $w(z)$ is inherited the same property.
Therefore by the definition
of the negation for absolute continuity
\cite{28} a constant $k$ can be chosen such that for any $\alpha>0$
there exists a positive constant $c$
and non overlapping intervals $(z'_i,z''_i)\in [-c,c]$, such that
$\sum_i|w(z''_i)-w(z'_i)|\ge k\quad\text{and}\quad \sum_i\int_{z'_i}^
{z''_i}n(z)dz\le\alpha.$
Put
$$
v^\alpha(z)=
\begin{cases}
-{1\over\sqrt{\alpha}}\text{sign} \ [w(z''_i)-w(z'_i)], & z'_i<z\le z''_i
\\
0, & \text{otherwise}.
\end{cases}
$$
Show that (8.3) holds. Evaluate from above $I(v^\alpha)$:
$$
\begin{aligned}
I(v^\alpha)
&
=\ell\int_\mathbb{R}\big[(v^\alpha(z))^2+|v^\alpha(z)|\big]n(z)dz+\int_\mathbb{R}w(z)
dv^\alpha(z)
\\
&
\le\ell\Big({1\over\alpha}+{1\over\sqrt{\alpha}}\Big)
\sum_i\int_{z'_i}^{z''_i}n(z)dz
+\sum_iw(z'_i)[v^\alpha(z''_i)-v^\alpha(z'_i)].
\\
&
\le \ell\big(1+\sqrt{\alpha}\big)
+\sum_iw(z'_i)[v^\alpha(z''_i)-v^\alpha(z'_i)].
\end{aligned}
$$
Now, summing by parts, we find
$\sum_iw(z'_i)[v^\alpha(z''_i)-v^\alpha(z'_i)]=
-\sum_iv^\alpha(z''_i)[w(z''_i)-w(z'_i)].$
On the other hand, from the definition of $v^\alpha(z)$ it follows
$\sum_iv^\alpha(z''_i)[w(z''_i)-w(z'_i)]=
{1\over\sqrt{\alpha}}\sum_i|w(z''_i)-w(z'_i)|\ge {k\over\sqrt{\alpha}}.$
Thereby
$$
I(v^\alpha)\le \ell(1+\sqrt{\alpha})-{k\over\sqrt{\alpha}}\to \infty, \
\alpha\to 0.
$$
Evidently to satisfy (8.4), it is sufficient to choose approximating functions
$v^\alpha_m(z), m$ $\ge 1$ which are compactly supported and continuously
differentiable and
such that $\lim_mv^\alpha_m(z)=v^\alpha(z)$ in every point of continuity
of $v^\alpha(z)$.

\medskip
\noindent
Assume $\nu\not\ll\lambda$. Put $K^\nu(dz)=
\int_0^TK_\nu(t,dz)dt$ and note that $K^\nu(dz)\not\ll dz.$
Use Lebesgue's decomposition:
$K^\nu(dz)=q(z)dz+K^\perp(dz)$, where $q(z)$ is a density of
absolutely continuous part of $K^\nu(dz)$ and $K^\perp(dz)$ is its
singular part. Taking $u(t,z)\equiv u(z)$ which is compactly supported, say,
on $[-c,c]$ we find
$$
\int_0^T\mathcal {D}u(z)\nu(dt,dz)=\int_{-c}^c\mathcal {D}u(z)q(z)dz+
\int_{-c}^c\mathcal {D}u(z)K^\perp(dz).
$$
Since $|u'(z)|\le |u'(0)|+\int_{-c}^c|u''(y)|dy$ there exists constant,
say, $\ell$, such that
\newline
$\int_{-c}^c\mathcal {D}u(z)q(z)dz$ $\le \ell(1+
\int_{-c}^c|u''(y)|dy)$ and so we arrive at an upper estimate
$$
\int_0^T\mathcal {D}u(z)\nu(dt,dz)\le \ell\Big(1+
\int_{-c}^c|u''(z)|dz\Big)+{1\over 2}\int_{-c}^c\sigma^2(z)u''(z)K^\perp(dz).
$$
Then, using the singularity of  $K^\perp(dz)$ and $dz$, one can choose $u''(z)$
such that the second integral is less any negative quantity while the first
remains bounded.

\medskip
\noindent
{\bf 3.} $\underline{\text{Approximation of rate function.}}$
\newline For `$dX=\dot{X}_tdt, \ d\nu=nd\lambda, \ d_zn=n'_zdz$'
denote by
$$
\begin{aligned}
&
S_T(X,\nu)=\int_0^T{[\dot{X}_t-A_\nu(t,X_t)]^2\over B^2_\nu(t,X_t)}dt,
\\
&
F_T(\nu)=\int_0^T\int_\mathbb{R}\sigma^2(z)\Big[{n'_z(t,z)\over n(t,z)}-
{p'(z)\over p(z)}\Big]^2n(t,z)dzdt.
\end{aligned}
$$
Also note one to one correspondence between density $n(t,z)$ and function
$v_\nu(t,z)$ defined in (6.1):
$$
n(t,z)=n(t,0){p(z)\over p(0)}\exp\Big(2\int_0^z{v_\nu(t,y)\over\sigma^2(y)}
dy\Big).
\eqno (8.5)
$$
Put
$$
\phi(t)=\int_\mathbb{R}|n'_z(t,y)|dy.
\eqno (8.6)
$$

\smallskip
{\bf Lemma A.4.} {\it Let $B^2(x,z)\ge \beta^2>0.$ If
$S_T(X,\nu)<\infty, \ F_T(\nu)<\infty$, then $\nu$ can be approximated by
a sequence of measures $\nu^{(k)}, k\ge 1$, satisfying the property:
$d\nu^{(k)}=n^{(k)}d\lambda, \ d_zn^{(k)}=n^{(k)}_zdz$, such that
the function
$v_{\nu^{(k)}}(t,z)$, corresponding to $n^{(k)}(t,z)$,
is compactly supported in $z$ and continuously differentiable in ($t,z$)
and what is more
$$
\begin{aligned}
&
\lim_k\rho_T(\nu,\nu^{(k)})=0
\\
&
\lim_kS_T(X,\nu^{(k)})=S_T(X,\nu)
\\
&
\lim_kF_T(\nu^{(k)})=F_T(\nu).
\end{aligned}
\eqno(8.7)
$$}

\medskip
\noindent
{\it Proof.} Introduce a chain of expanding subclasses of measures $\nu$
characterized in terms of $n(t,z)$ and $v_\nu(t,z)$:

0) $v_\nu(t,z)$ is compactly supported in $z$ and continuously differentiable
in ($t,z$);

1) $v_\nu(t,z)$ is compactly supported in $z$ and bounded;

2) $v_\nu(t,z)$ is compactly supported,
$\inf_{t\le T, z\in \mathbb{R}}{n(t,z)\over p(z)}>0$
and $\sup_{t\le T}[n(t,0)+ \phi(t)]<\infty$;

3) $v_\nu(t,z)$ is compactly supported,
$\inf_{t\le T, z\in \mathbb{R}}{n(t,z)\over p(z)}>0$;

4) $v_\nu(t,z)$ is compactly supported;

5) $v_\nu(t,z)$ satisfies the assumptions of the lemma.

\medskip
\noindent
The proof is based on the following fact. If measure $\nu$ from calss `i'
(i=1,...,5) can be approximated by $\nu^{(k)}, k\ge 1$ from class `i-1'
in a sense (8.7), then the statement of the lemma holds.

\medskip
\noindent
Assume $\nu^{(k)}, k\ge 1$ is such that
$$
\begin{aligned}
&
\Lambda_T-\lim_kn^{(k)}(t,z)=n(t,z), \ (\Lambda_T(dt,dz)=I_{[0,T]}dtdz),
\\
&
\lim_kF_T(\nu^{(k)})=F_T(\nu).
\end{aligned}
\eqno(8.8)
$$
Then, by Scheffe's theorem, we have
$\lim_k\int_0^T\int_\mathbb{R}|n(t,z)-n^{(k)}(t,z)|dtdz=0$
that is $\nu^{(k)}$ converges to $\nu$ in the total variation norm which
implies convergence in Levy-Prohorov's metric too:
$\rho_T(\nu,\nu^{(k)})\to 0$. Since
$$A_{\nu^{(k)}}(t,X_t)=
\int_\mathbb{R}A(X_t,z)n^{(k)}(t,z)dz, \
B^2_{\nu^{(k)}}(t,X_t)=
\int_\mathbb{R}B^2(X_t,z)n^{(k)}(t,z)dz
$$
by Lebesque dominated theorem
$S_T(X,\nu^{(k)})\to S_T(X,\nu)$. Therefore, for all steps
of approximations only (8.8) has to be checked.

\medskip
\noindent
Assume $\nu$ is from class `1'.
Approximate $v_\nu(t,z)$ by $v_\nu^{(k)}(t,z)$:
$$
\lim_k\int_0^T\int_\mathbb{R}[v_\nu(t,z)-v_\nu^{(k)}(t,z)]^2
\big(1+n(t,z)\big)dtdz=0,
$$
where for all $k$ the functions $v_\nu^{(k)}(t,z), k\ge 1$  are
compactly supported continuously differentiable in $(t,z)$.
Without loss of a generality one can assume that all function
are bounded by the same constant.
Similarly to (8.5) define a density of $\nu^{(k)}$:
$$
n^{(k)}(t,z)=n^{(k)}(t,0){p(z)\over p(0)}\exp\Big(\int_0^z{v^{(k)}_\nu(t,y)
\over\sigma^2(z)}dy\Big),
\eqno (8.9)
$$
with
$
n^{(k)}(t,0)=\Big(\int_\mathbb{R}{p(z)\over p(0)}\exp\Big(\int_0^z{v^{(k)}_\nu(t,y)\over
\sigma^2(z)}dy\Big)dz\Big)^{-1}.
$
Put
$$
\nu^{(k)}(dt,dz)=n^{(k)}(t,z)dtdz.
$$
Evidently $\nu^{(k)}$ belongs to class`0'.
It is easy to check that
$$
v^{(k)}_\nu(t,z)\equiv v_{\nu^{(k)}}(t,z)
\eqno (8.10)
$$
and the validity of the first part in (8.8). To verify the second part in
(8.8), note that
$F_T(\nu^{(k)})=4\int_0^T\int_\mathbb{R}{(v_{\nu^{(k)}}(t,z))^2
\over \sigma^2(z)}n^{(k)}(t,z)dzdt$ and consequently
$$
\begin{aligned}
\big|F_T(\nu)-F_T(\nu^{(k)})\big|
&
\le \text{const.}\int_0^T\int_\mathbb{R}|n(t,z)-n^{(k)}(t,z)|dtdz
\\
&
\quad +\text{const.}\int_0^T\int_\mathbb{R}\big|v^2_\nu(t,z)-(v_{\nu^{(k)}}(t,z))^2
\big|dtdz
\\
&
\to 0, \ k\to\infty.
\end{aligned}
$$

\medskip
\noindent
Assume $\nu$ is from class `2'. For the definiteness assume that
there exists positive constant $z_\circ$ such that $v_\nu(t,z)\equiv 0$
out of $[-z_\circ,z_\circ]$.
Put $v_\nu^{(k)}(t,z)=v_\nu(t,z)I(|n'_z(t,z)|\le k)$, define
$n^{(k)}(t,z)$ by (8.9) and take $\nu^{(k)}$ with this density. It belongs to
class `1'.
Herewith, $v_\nu^{(k)}(t,z)$ is defined by (8.10). It is clear that
the first part in (8.8) holds and below
we check the validity of the second part. We have
$$
\begin{aligned}
&
F_T(\nu^{(k)})
=4\int_0^T\int_{|z|\le z_\circ}{v_\nu^2(t,z)\over \sigma^2(z)}
I(|n_z(t,z)|\le k)n^{(k)}(t,z)dzdt
\\
&
F_T(\nu)
=4\int_0^T\int_{|z|\le z_\circ}{v_\nu^2(t,z)\over \sigma^2(z)}
n(t,z)dzdt.
\end{aligned}
$$
The required convergence $F_T(\nu^{(k)})\to F_T(\nu)$ holds by Lebesgue
dominated theorem since
$
n^{(k)}(t,z)\le p(z)\exp\big(2\phi(t)\big)\le
\text{const.}n(t,z).
$

\medskip
\noindent
Assume $\nu$ is from class `3'. Putting $v_\nu^{(k)}(t,z)=v_\nu(t,z)
I(n(t,0)+\phi(t)\le k)$ we arrive at
$$
n^{(k)}(t,z)=
\begin{cases} n(t,z), & n(t,0)+\phi(t)\le k\\
p(z), & n(t,0)+\phi(t)>k.
\end{cases}
$$
and since  $n^{(k)}(t,0)\le k+p(0)$ and $\phi^{(k)}(t)\le k+\int_\mathbb{R}|p'(z)|dz$
measure $\nu^{(k)}$ with density $n^{(k)}(t,z)$ belongs to class `2'.
It is clear that the first part in (8.8) holds and
$$
\begin{aligned}
|F_T(\nu)-F_T(\nu^{(k)})|
&
=4\int_0^T\int_{|z|\le z_\circ}{v_\nu^2(t,z)
\over \sigma^2(z)}I(n(t,0)+\phi(t)>k)p(z)dzdt
\\
&
\le\text{const.}\int_0^T\int_{|z|\le z_\circ}{v_\nu^2(t,z)
\over \sigma^2(z)}I(n(t,0)+\phi(t)>k)n(t,z)dzdt
\\
&
\to 0, \ k\to\infty.
\end{aligned}
$$

\medskip
\noindent
Assume $\nu$ is from class `4'. Put
$n^{(k)}(t,z)=c^{(k)}(t)\big(n(t,z)\vee p(z)\big),$
where $c^{(k)}(t)=\Big(\int_\mathbb{R}\big(n(t,z)\vee p(z)\big)dz\Big)^{-1}$ is
norming constant. $\nu^{(k)}$ with this density belongs to class `3'.
The first part in (8.8)
holds and what is more $\lim_kc^{(k)}(t)=1$. On the other hand, since
$v_{\nu^{(k)}}(t,z)=v_\nu(t,z)I(n(t,z)\ge p(z)/k)$ we obtain
$$
\begin{aligned}
F_T(\nu^{(k)})
&
=4\int_0^T\int_\mathbb{R}{v_\nu^2(t,z)\over \sigma^2(z)}
I(n(t,z)\ge p(z)/k)c^{(k)}(t)n(t,z)dzdt
\\
&
\to 4\int_0^T\int_\mathbb{R}
{v_\nu^2(t,z)\over \sigma^2(z)}n(t,z)dzdt=F_T(\nu).
\end{aligned}
$$

\medskip
\noindent
Assume $\nu$ is from class `5'. Put
$v_\nu^{(k)}(t,z)=v_\nu(t,z)T(|z|\le k)$ and define $n^{(k)}(t,z)$ by
(8.9). Then
$$
n^{(k)}(t,z)=n^{(k)}(t,0)
\begin{cases}
p(z){n(t,k)\over p(k)}, & z>k\\
n(t,z), & |z|\le k\\
p(z){n(t,-k)\over p(-k)}, & z<-k.
\end{cases}
$$
Taking $\nu^{(k)}$ with this density and noticing that
$\lim_kn^{(k)}(t,0)=1$
we find
$$
F_T(\nu^{(k)})=4\int_0^T\int_{|z|\le k}
{\big(v_\nu(t,z)\big)^2\over \sigma^2(z)}
c^{(k)}(t)n(t,z)dzdt\to F_T(\nu),
$$
i.e. both parts in (8.8) hold.

\medskip
\noindent\medskip
\noindent
{\bf 4.} $\underline{\text{Ergodic property.}}$
\newline
Consider diffusion pair $(\widetilde{X}_t^\varepsilon,
\widetilde{\xi}_t^\varepsilon)$ defined by It\^o's differential equations
w.r.t. independent Wiener processes $W_t$ and $V_t$:
$$
\begin{aligned}
&
d\widetilde{X}_t^\varepsilon=G(t,\widetilde{X}_t^\varepsilon,
\widetilde{\xi}_t^\varepsilon)dt+\sqrt{\varepsilon}
B(\widetilde{X}_t^\varepsilon,\widetilde{\xi}_t^\varepsilon)dW_t
\\
&
d\widetilde{\xi}_t^\varepsilon={1\over \varepsilon}b(t,
\widetilde{\xi}_t^\varepsilon)dt+
{1\over\sqrt{\varepsilon}}\sigma(\widetilde{\xi}_t^\varepsilon)dV_t
\end{aligned}
\eqno(8.11)
$$
subject to $(x_0,z_0)$, where $B(x,z)$ and $\sigma(z)$ are
functions involving in (1.1).  Assume $b(t,z)$ is
continuous it $(t,z)$, continuously differentiable in $t$,
Lipschitz continuous in $z$ uniformly in t, and $zb(t,z)$ is
negative for large $|z|$ uniformly in $t$. Also assume that
$$
G(t,x,z)={\dot{X}_t-A_p(t,x)\over B_p(t,x)}B(x,z)+A(x,z),
\eqno (8.12)
$$
where $A(x,z)$ involves in (1.1), $\dot{X}_t$ is the Radon-Nykodim
derivative of absolute continuous function $X_t$ from $\Bbb{C}$ with
$X_0=x_0$, $A_p(t,x)=\int_\mathbb{R}A(x,z)p(t,z)dz$,
$B_p(t,x)=\sqrt{\int_\mathbb{R}B^2(x,z)p(t,z)dz}$, and where (comp. (2.1))
$$
p(t,z)=c(t){\exp\Big(2\int_0^z{b(t,y)\over \sigma^2(y)}dy\Big)\over
\sigma^2(z)}
$$
with norming function $c(t)$ such that $\int_\mathbb{R}p(t,z)dz=1.$ Introduce an
occupation
measure $\widetilde{\nu}^\varepsilon(dt,dz)$:
$\widetilde{\nu}^\varepsilon(\Delta\times\Gamma)=
\int_0^\infty I(t\in \Delta, \widetilde{\xi}_t^\varepsilon\in\Gamma)dt$
and put $\nu(dt,dz)=p(t,z)dzdt$.

\smallskip{\bf Lemma A.5.} {\it
$$
\mathsf{P}-\lim_{\varepsilon\to 0}\rho_T(
\widetilde{\nu}^\varepsilon,\nu)=0
\quad\text{and}\quad
\mathsf{P}-\lim_{\varepsilon\to 0}r_T(
\widetilde{X}^\varepsilon,X)=0
$$
}
\medskip
\noindent
{\it Proof.} It is clear, the first statement of the lemma
is equivalent to: for any bounded and continuous function $h(t,z)$
$\int_0^T\int_\mathbb{R}h(t,
\widetilde{\xi}_t^\varepsilon)dt\to \int_0^T\int_\mathbb{R}h(t,z)p(t,z)dzdt$
in probability or, for $h^\circ(t,z)=h(t,z)-\int_0^T\int_\mathbb{R}h(t,y)p(t,y)dydt$,
$$
{\mathsf P}-\lim_{\varepsilon\to 0}\int_0^T\int_\mathbb{R}h^\circ(t,
\widetilde{\xi}_t^\varepsilon)dt=0.
$$
First we check it for continuously differentiable
in $t,z$ function $h(t,z)$, having bounded partial derivatives.
Straightforward calculation brings Kolmogorov's type
differential equation ($t$ is fixed):
$$
{1\over 2}{\partial \over \partial z}(\sigma^2(x)p(t,x))=
b(t,z)p(t,z).
$$
A conjugate differential equation
$$
{1\over 2}\sigma^2(z){\partial w(t,z)\over \partial z}+b(t,z)w(t,z)=
h^\circ(t,z).
\eqno (8.13)
$$
obeys a solution
$$
w(t,z)={2\over \sigma^2(z)p(t,z)}\int_{-\infty}^zh^\circ(t,y)p(t,y)dy.
$$
It is clear that properties of $h(t,z)$ are
inherited by $w(t,z)$ and so function $u(t,z)=\int_0^zw(t,y)dy$ is
continuously differentiable once in $t$
and twice in $z$ and what is more, due to the boundedness of $w(t,z)$,
there exists a positive constant, say $\ell$, such that
$|u(t,z)|\le \ell |z|\quad \text{and}\quad |u_t(t,z)|\le \ell|z|.$
Applying It\^o's formula to $u(t,\widetilde{\xi}_t
^\varepsilon)$ and taking into account that $w(t,z)$ is solution of
differential equation (8.13) we find
$u(T,\widetilde{\xi}^\varepsilon_T)
=u(0,\xi_0)+\int_0^Tu'_t(t,\widetilde{\xi}^\varepsilon_t)dt$
\newline
$+{1\over\sqrt{\varepsilon}}\int_0^Tw(t,\widetilde{\xi}^\varepsilon_t)
\sigma(\widetilde{\xi}^\varepsilon_t)dV_t+
{1\over\varepsilon}\int_0^Th^\circ(t,\widetilde{\xi}^\varepsilon_t)dt$
that is
$$
\begin{aligned}
\int_0^Th^\circ(t,\widetilde{\xi}^\varepsilon_t)dt
&
=\varepsilon u(T,\widetilde{\xi}^\varepsilon_T)
-\varepsilon u(0,\xi_0)
\\
&
\quad-\varepsilon\int_0^Tu'_t(t,\widetilde{\xi}^\varepsilon_t)dt
-\sqrt{\varepsilon}\int_0^Tw(t,\widetilde{\xi}^\varepsilon_t)
\sigma(\widetilde{\xi}^\varepsilon_t)dV_t.
\end{aligned}
\eqno(8.14)
$$
The second term in the right hand of (8.10) converges to zero;
the last term converges to zero in probability since by Problem 1.9.2
in \cite{23} the mentioned convergence follows from
$\varepsilon\int_0^Tw^2(t,\widetilde{\xi}^\varepsilon_t)
\sigma^2(\widetilde{\xi}^\varepsilon_t)dt\to 0$; other two terms  converge
to zero in probability if
$\lim_{\varepsilon\to 0}\varepsilon^2{\mathsf E}\sup_{t\le T}(
\widetilde{\xi}_t^\varepsilon)^2=0.$
To check the last, apply It\^o's formula to
$(\varepsilon\widetilde{\xi}_t^\varepsilon)^2$:
$$
\begin{aligned}
(\varepsilon\widetilde{\xi}_t^\varepsilon)^2
&
=(\varepsilon\xi_0)^2+
2\varepsilon\int_0^tb(s,\widetilde{\xi}_s^\varepsilon)
\widetilde{\xi}_s^\varepsilon ds+
\varepsilon\int_0^t\sigma^2(\widetilde{\xi}_s^\varepsilon)ds
&
\quad+2\varepsilon^{3/2}\int_0^t\widetilde{\xi}_s^\varepsilon
\sigma(\widetilde{\xi}_s^\varepsilon)dV_s.
\end{aligned}
$$
The function $b(s,z)$ is such that $zb(t,z)$ is negative for large $|z|$
what implies
$(\varepsilon\widetilde{\xi}_t^\varepsilon)^2\le (\varepsilon\xi_0)^2+
T\varepsilon\text{const.}
+2\varepsilon^{3/2}\int_0^t\widetilde{\xi}_s^\varepsilon
\sigma(\widetilde{\xi}_s^\varepsilon)dV_s$. Thereby
$
{\mathsf E}(\varepsilon\widetilde{\xi}_t^\varepsilon)^2\le (\varepsilon\xi_0)^2+
T\varepsilon\text{const.}$ In turn, using Doob's inequality
(see e.g. Theorem 1.9.1 in \cite{23}), we arrive at
${\mathsf E}\sup_{t\le T}(\varepsilon\widetilde{\xi}_t^\varepsilon)^2
\le (\varepsilon\xi_0)^2+T\varepsilon\text{const.}+
\text{const.}\varepsilon^3\int_0^T
{\mathsf E}(\varepsilon\widetilde{\xi}_t^\varepsilon)^2dt$ and, due to
the obtained above upper bound for
${\mathsf E}(\varepsilon\widetilde{\xi}_t^\varepsilon)^2$, the result holds.

\medskip
\noindent
If $h(t,z)$ is
bounded and continuous only, it can be approximated by smooth functions
$h_m(t,z), m\ge 1$ in the following sense: for any $k\ge 1$
$$
\lim_m\sup_{t\le T, |z|\le k}|h(t,z)-h_m(t,z)|=0.
$$
Since for every $h_m(t,z)$ the statement of the lemma is proved, it
holds for $h(t,z)$ if
$$
\begin{aligned}
&
\lim_m\int_0^T\int_\mathbb{R}|h(t,z)-h_m(t,z)|p(t,z)dzdt=0
\\
&
{\mathsf P}-\lim_m\bar{\lim_{\varepsilon\to 0}}
\int_0^T|h(t,\widetilde{\xi}_t^\varepsilon)-
h_m(t,\widetilde{\xi}_t^\varepsilon)|n(t,z)dt=0.
\end{aligned}
$$
The first takes place since $\lim_k\int_0^T\int_{|z|>k}n(t,z)dzdt=0$
while the second from
${\mathsf P}-\lim_m\bar{\lim}_{\varepsilon\to 0}
\int_0^TI(|\widetilde{\xi}_t^\varepsilon)|>k)dt=0$
and the fact that one can choose smooth bounded
functions $g_k(z), k\ge 1$ such that $I(|z|>k)\le g_k(z)$, $\lim_kg_k(z)=0,
z\in \mathbb{R}$ and by proved above
$\int_0^Tg_k(\widetilde{\xi}_t^\varepsilon)dt\to\int_0^T\int_\mathbb{R}g_k(z)p(t,z)dzdt
\to 0, \ k\to\infty.$
\medskip
\noindent
To check the second statement, put $\Delta_t=\widetilde{X}_t^\varepsilon-X_t$.
From the first equation in (8.11) we find
$$
\begin{aligned}
\Delta_t
&
=\int_0^t\Big[{\dot{X}_t-A_p(s,X_s)\over B_p(s,X_s)}
B(\widetilde{X}_s^\varepsilon,\widetilde{\xi}_s^\varepsilon)+
A(\widetilde{X}_s^\varepsilon,\widetilde{\xi}_s^\varepsilon)-\dot{X}_s
\Big]ds
\\
&
\quad+\sqrt{\varepsilon}\int_0^t
B(\widetilde{X}_s^\varepsilon,\widetilde{\xi}_s^\varepsilon)dW_s
\\
&
=\int_0^t\Big[{\dot{X}_t-A_p(s,X_s)\over B_p(s,X_s)}
\big(B(\widetilde{X}_s^\varepsilon,\widetilde{\xi}_s^\varepsilon)-
B(X_s,\widetilde{\xi}_s^\varepsilon)\big)\Big]ds
\\
&
\quad+\int_0^t\Big[{\dot{X}_s-A_p(s,X_s)\over B_p(s,X_s)}
\big(B(X_s,\widetilde{\xi}_s^\varepsilon)-
B_p(s,X_s)\big)\Big]ds
\\
&
\quad+\int_0^t\big(A(\widetilde{X}_s^\varepsilon,
\widetilde{\xi}_s^\varepsilon)-
A(X_s,\widetilde{\xi}_s^\varepsilon)\big)ds
\\
&
\quad+\int_0^t\big(A(X_s,\widetilde{\xi}_s^\varepsilon)-
A_p(s,X_s)\big)ds
\\
&
\quad+\sqrt{\varepsilon}\int_0^t
B(\widetilde{X}_s^\varepsilon,\widetilde{\xi}_s^\varepsilon)dW_s.
\end{aligned}
$$
For brevity put $\varphi_s=
{\dot{X}_s-A_p(s,X_s)\over B_p(s,X_s)}$. Then by the Lipschitz
continuity of functions $A(x,z),B(x,z)$ in $x$
uniformly in $z$, say, with constant $\ell$, we obtain
$$
\begin{aligned}
|\Delta_t|
&
\le\ell \int_0^t(1+|\varphi_s|)\Delta_sds
\\
&
\quad+\sup_{t\le T}\Big|\int_0^t\varphi_s\big(B(X_s,\widetilde{\xi}_s^\varepsilon)-
B_p(s,X_s)\big)ds\Big|
\\
&
\quad+\sup_{t\le T}\Big|\int_0^t\big(A(X_s,\widetilde{\xi}_s^\varepsilon)-
A_p(s,X_s)\big)ds\Big|
\\
&
\quad+\sqrt{\varepsilon}\sup_{t\le T}\Big|\int_0^t
B(\widetilde{X}_s^\varepsilon,\widetilde{\xi}_s^\varepsilon)dW_s\Big|.
\end{aligned}
$$
Therefore, by Gronwall-Bellman's inequality
$$
\begin{aligned}
\sup_{t\le T}|\Delta_t|
&
\le \exp\Big(\ell \int_0^T(1+|\varphi_s|)ds\Big)
\\
&
\quad\times\Big[\sup_{t\le T}\Big|\int_0^t\varphi_s\big(B(X_s,
\widetilde{\xi}_s^\varepsilon)-B_p(s,X_s)\big)ds\Big|
\\
&
\quad+\sup_{t\le T}\Big|\int_0^t\big(A(X_s,\widetilde{\xi}_s^\varepsilon)-
A_p(s,X_s)\big)ds\Big|
\\
&
\quad+\sqrt{\varepsilon}\sup_{t\le T}\Big|\int_0^t
B(\widetilde{X}_s^\varepsilon,\widetilde{\xi}_s^\varepsilon)dW_s\Big|\Big]
\end{aligned}
$$
Hence, the second statement holds if
$$
{\mathsf P}-\lim_{\varepsilon\to 0}\sup_{t\le T}
\sqrt{\varepsilon}\sup_{t\le T}\Big|\int_0^t
B(\widetilde{X}_s^\varepsilon,\widetilde{\xi}_s^\varepsilon)dW_s\Big|=0
\eqno (8.15)
$$
and for any measurable
function $\psi_s$ such that $\int_0^T\psi_s^2ds<\infty$ and any continuous
function $C(x,z)$, being Lipschitz's continuous in $x$ uniformly in $z$,
$$
{\mathsf P}-\lim_{\varepsilon\to 0}\sup_{t\le T}\Big|\int_0^t\psi_s\big(
C(X_s,\widetilde{\xi}_s^\varepsilon)-
C_p(s,X_s)\big)ds\Big|=0,
\eqno (8.16)
$$
where $C_p(s,X_s)=\int_\mathbb{R}C(X_s,z)p(s,z)dz.$ It can be shown (see e.g.
the method of proving the statement (2) of
Theorem 4.6 Ch. 4 in \cite{24}) that $\sup_{t\le T}{\mathsf E}(\widetilde{X}_t
^\varepsilon)^2\le \text{const.}$ and so
${\mathsf E}\int_0^TB^2(\widetilde{X}_s^\varepsilon,
\widetilde{\xi}_s^\varepsilon)ds\le \text{const.}$ Consequently, by Doob's
inequality (see e.g. Theorem 1.9.1 \cite{23}) we get
${\mathsf E}\sup_{t\le T}\Big|
\sqrt{\varepsilon}\int_0^t
B(\widetilde{X}_s^\varepsilon,\widetilde{\xi}_s^\varepsilon)dW_s\Big|^2
\le \varepsilon\text{const.}$ that is  (8.15) holds.
To check the validity of (8.16) with $\psi_sC(x,z)\ge 0$, note that, due to
Problem 5.3.2 in\cite{23}, it is sufficient to prove
$$
{\mathsf P}-\lim_{\varepsilon\to 0}\int_0^t\psi_s\big(
C(X_s,\widetilde{\xi}_s^\varepsilon)-
C_p(s,X_s)\big)ds=0, \ \forall t\le T
\eqno (8.17)
$$
and what is more, due to an arbitrariness of $\psi_s$ and $C(x,z)$, (8.17)
implies (8.16) in the general case since one can use separately (8.17)
for positive $(\psi_sC(x,z))^+$ and negative $(\psi_sC(x,z))^-$ parts.
Therefore, (8.17) remains to be verified. If $\psi_s$ is bounded
and continuous, (8.17) takes place by virtue of the first statement
of the lemma. If  only
$\int_0^T\psi_s^2ds<\infty$, approximate $\psi_s$ by bounded and
continuous functions $\psi_s^{(k)}, \ k\ge 1$ such that
$\lim_k\int_0^T(\psi_s-\psi_s^{(k)})^2ds=0$ and, due to the boundedness
in $z$ of $C(x,z)$ and Cauchy-Schwartz's inequality, we find that
$$
\Big|\int_0^t(\psi_s-\psi_s^{(k)})\big(
C(X_s,\widetilde{\xi}_s^\varepsilon)-
C_p(s,X_s)\big)ds\Big|
$$
$$\le\text{const.}\sqrt{\int_0^T(\psi_s-\psi_s^{(k)})^2ds
}\to 0, \ k\to\infty
$$
that is (8.17) takes place since it holds for every $\psi_s^{(k)}.$

\medskip
\noindent
{\bf 5.} $\underline{\text{LD-regularization.}}$
\newline
Parallel to $X_t^\varepsilon$, defined in (1.1), determine new diffusion
$X_t^{\varepsilon,\beta}$ with uniformly non singular diffusion parameter
$B^2(x,z)+\beta^2, \beta^2>0$, letting
$$
dX_t^{\varepsilon,\beta}=A(X_t^{\varepsilon,\beta},\xi^\varepsilon_t)dt+
\sqrt{\varepsilon}\big[B(X_t^{\varepsilon,\beta},\xi^\varepsilon_t)dW_t+
\beta dW'_t\big]
\eqno (8.18)
$$
subject to the same initial point $x_0$, where $W'_t$ is a Wiener process
independent of $(W_t, \xi_t^\varepsilon).$

\smallskip
{\bf Lemma A.6.} {\it Under assumption (A.1) for every $T>0$ and $\eta>0$
$$
\lim_{\beta\to 0}\bar{\lim}_{\varepsilon\to 0}\varepsilon\log\mathsf{P}(
r_T(X^{\varepsilon,\beta},X^\varepsilon)>\eta)=-\infty.
$$
}

\medskip
{\it Proof.} Put $\Delta_t=X_t^{\varepsilon,\beta}-X_t^\varepsilon$,
and
$$
a_1(x',x'',z)={A(x'',z)-A(x',z)\over x''-x'},\quad
a_2(x',x'',z)={B(x'',z)-B(x',z)\over x''-x'},
$$
where for $x'=x''$ $a_i(x',x',z), i=1,2$
are Radon-Nikodym's derivatives. Note that for $x'\neq x''$ $a_i(x',x'',z),$
are bounded, say, by constant $\ell$,  and so $a_i(x',x',z)$
can be taken bounded by the same constant. For brevity, denote by
$\alpha_i(t)=a_i(X_t^{\varepsilon,\beta},X_t^\varepsilon,\xi_t^\varepsilon),
i=1,2$. (8.18) and (1.1) imply:
$$
\Delta_t=\int_0^t\alpha_1(s)\Delta_sds+\sqrt{
\varepsilon}\int_0^t\alpha_2(s)\Delta_sdW_s+\sqrt{\varepsilon}\beta W'_t.
$$
Letting $\mathcal {E}_t=\exp\Big(\int_0^t[\alpha_1(s)-(\varepsilon/2)
\alpha^2_2(s)]ds+\sqrt{\varepsilon}\int_0^t\alpha_2(s)dW_s\Big)$ and
using It\^o's formula, we find that $\Delta_t=\sqrt{\varepsilon}\beta\mathcal {E}_t
\int_0^t\mathcal {E}_s^{-1}dW'_s$ and thereby
$$
\sup_{t\le T}\big|\Delta_t\big|\le \sqrt{\varepsilon}\beta\sup_{t\le T}\mathcal {E}_t
\sup_{t\le T}\Big|\int_0^t\mathcal {E}_s^{-1}dW'_s\Big|.
$$
Put $\Gamma_N=\{{1\over N}\le \inf_{t\le T}\mathcal {E}_t\le \sup_{t\le T}\mathcal {E}_t
\le N\}$ and use an upper estimate
$$
\begin{aligned}
\mathsf{P}(\sup_{t\le T}\big|\Delta_t\big|>\eta)
&
\le \mathsf{P}(\sup_{t\le T}\big|\Delta_t\big|>\eta, \Gamma_N)+
\mathsf{P}(\Omega\setminus\Gamma_N)
\\
&
\le 2\max\Big[\mathsf{P}(\sup_{t\le T}\big|\Delta_t\big|>\eta, \Gamma_N),
\mathsf{P}(\Omega\setminus\Gamma_N)\Big]
\end{aligned}
$$
which implies, due to the boundedness of $\alpha_i(s), i=1,2$,
the desired statement if
$$
\begin{aligned}
&
\lim_N\bar{\lim}_{\varepsilon\to 0}\varepsilon\log\mathsf{P}
(\sqrt{\varepsilon}\sup_{t\le T}\Big|\int_0^t\alpha_2(s)dW_s\Big|>N)=-\infty
\\
&
\lim_{\beta\to 0}\bar{\lim}_{\varepsilon\to 0}\varepsilon\log\mathsf{P}
(\sqrt{\varepsilon}\beta\sup_{t\le T}\Big|\int_0^t\mathcal {E}_s^{-1}
dW'_s\Big|>\eta, \Gamma_N)=-\infty, \forall N\ge 1.
\end{aligned}
\eqno(8.19)
$$
Let $\tau=\{t\le T:
\Big|\int_0^t\alpha_2(s)dW_s\Big|>(N/\sqrt{\varepsilon})\}$ and
$\sigma=\{t\le t:\Big|\int_0^t\mathcal {E}_s^{-1}
dW'_s\Big|>(\eta/\sqrt{\varepsilon}\beta)\}$. It is clear
that (8.19) is equivalent to:
$$
\begin{aligned}
&
\lim_N\bar{\lim}_{\varepsilon\to 0}\varepsilon\log\mathsf{P}
\big(\sqrt{\varepsilon}\int_0^\tau\alpha_2(s)dW_s\ge N
\ (\text{or} \ \le -N)\big)=-\infty
\\
&
\lim_{\beta\to 0}\bar{\lim}_{\varepsilon\to 0}\varepsilon\log\mathsf{P}
\big(\sqrt{\varepsilon}\beta\int_0^\sigma\mathcal {E}_s^{-1}
dW'_s\ge \eta \ (\text{or} \ \le -\eta), \Gamma_N\big)=-\infty, \forall N\ge 1.
\end{aligned}
\eqno(8.20)
$$
Below we check (8.20). To this end with $\lambda\in \mathbb{R}$,
introduce continuous local martingales:
$Z_t^1=\exp\Big(\lambda\int_0^t\alpha_2(s)dW_s-(\lambda^2/2)
\int_0^t\alpha_2^2(s)ds\Big)$ and $Z_t^2=\exp\Big(\lambda\int_0^t\mathcal {E}_s^{-1}
dW'_s-(\lambda^2/2)\int_0^t\mathcal {E}_s^{-2}ds\Big)$, where each of them is a
supermartingale too (see Problem 1.4.4. in \cite{23}) that is $\mathsf{E}
Z^1_\tau\le 1$ and $\mathsf{E}Z^2_\sigma\le 1$. Then we use obvious
inequalities:
$$
\mathsf{E}I(\sqrt{\varepsilon}\int_0^\tau\alpha_2(s)dW_s\ge N)
Z^1_\tau\le 1, \ \mathsf{E}I(\sqrt{\varepsilon}\beta\int_0^\sigma\mathcal {E}_s
dW'_s\ge \eta, \Gamma_N)Z^2_\sigma\le 1.
$$
Since for $\lambda>0$ we have
$\log Z^1_\tau\ge (\lambda N/\sqrt{\varepsilon})-(\lambda^2\ell^2T/2)$
and $\log Z^1_\sigma\ge (\lambda \eta/\sqrt{\varepsilon}\beta)-
(\lambda^2N^2T/2)$ on sets $\{\sqrt{\varepsilon}\int_0^\tau\alpha_2(s)dW_s
\ge N\}$ and $\{\sqrt{\varepsilon}\beta\int_0^\sigma\mathcal {E}_s
dW'_s\ge \eta, \Gamma_N\}$ respectively, (8.20) holds for
`$\ge N$' and `$\ge\eta$', by taking $\lambda^1=(N/\sqrt{\varepsilon}\ell^2T)$
and $\lambda^2=(\eta/\sqrt{\varepsilon}\beta N^2T)$. For
`$\le -N$' and `$\le-\eta$' the validity of (8.20) is proved in the same way.

\medskip
\noindent
{\bf Acknowledgements :} The author gratefully acknowledge the careful review
by the anonymous referee, who also pointed out a mistake in the original
version. Due to his comments the paper could be significantly improved.

\enddocument